\input amstex
\documentstyle{amsppt}
\magnification=\magstep 1
\catcode`\@=11
\def\nologo{\let\logo@=\relax}
\catcode`\@=\active
\nologo
\vsize6.95in

\topmatter
\title Divisors over determinantal rings defined by two by two minors \endtitle 
  \leftheadtext{Andrew R. Kustin}
\rightheadtext{Two by Two}
 \author Andrew R. Kustin\endauthor
 \address
Mathematics Department,
University of South Carolina,
Columbia, SC 29208\endaddress
\email kustin\@math.sc.edu \endemail

\abstract  
  Let $E$ and $G$ be free modules of rank $e$ and $g$, respectively, over a commutative noetherian ring $R$. The identity map on $E^*\otimes G$ induces the Koszul complex $$ \to \operatorname{Sym}_mE^*\otimes \operatorname{Sym}_nG\otimes {\tsize \bigwedge}^p(E^*\otimes G)\to \operatorname{Sym}_{m+1}E^*\otimes \operatorname{Sym}_{n+1}G\otimes {\tsize \bigwedge}^{p-1}(E^*\otimes G)\to   $$ and its dual $$\eightpoint \dots \to D_{m+1}E\otimes D_{n+1}G^*\otimes {\tsize \bigwedge}^{p-1}(E\otimes G^*)\to D_mE\otimes D_nG^*\otimes {\tsize \bigwedge}^p(E\otimes G^*)\to \dots\ .$$Let $\operatorname{H}_{\Cal N}(m,n,p)$ be the homology of the top complex at $\operatorname{Sym}_mE^*\otimes \operatorname{Sym}_nG\otimes {\tsize \bigwedge}^p(E^*\otimes G)$ and $\operatorname{H}_{\Cal M}(m,n,p)$  the homology of the bottom complex at $D_mE\otimes D_nG^*\otimes {\tsize \bigwedge}^p(E\otimes G^*)$. It is known  that $\operatorname{H}_{\Cal N}(m,n,p)\cong \operatorname{H}_{\Cal M}(m',n',p')$,
provided $m+m'=g-1$, $n+n'=e-1$,  $p+p'=(e-1)(g-1)$, and
$1-e\le m-n\le g-1$. In this paper we exhibit an explicit quasi-isomorphism $M$ of complexes which gives rise to this isomorphism. The mapping cone of $M$ is a split exact complex. Our complexes may be formed over the ring of integers; they can be passed to an arbitrary ring or field by base change.   

Knowledge of the homology of the top complex is equivalent to knowledge of the modules in the resolution of the  Segre module  $\operatorname{Segre}(e,g,\ell)$, for $\ell=m-n$. The  modules $\{\operatorname{Segre} (e,g,\ell)|\ell\in \Bbb Z\}$
are a set of representatives of the divisor class group of the determinantal ring defined by the $2 \times 2$ minors of an $e \times g$ matrix of indeterminants. If $R$ is the ring of integers, then the  homology $\operatorname{H}_{\Cal N}(m,n,p)$ is not always a free abelian group. In other words, if $R$ is a field, then the dimension of $\operatorname{H}_{\Cal N}(m,n,p)$ depends on the characteristic of $R$. The module $\operatorname{H}_{\Cal N}(m,n,p)$ is known when $R$ is a field of characteristic zero; however, this module is not yet known over arbitrary fields. 

The homology $\operatorname{H}_{\Cal N}(m,n,p)$ is equal to the homology of the simplicial complex known as a chessboard complex with multiplicities. A chessboard complex is the matching complex of a a complete bipartite graph.

The modules in the minimal resolution of the universal ring for finite length modules of projective dimension two are equal to modules of the form $\operatorname{H}_{\Cal N}(m,n,p)$.
 \endabstract

\endtopmatter

\document

\null\footnote""{2000 {\it Mathematics Subject Classification.} 13D25.}\footnote""{{\it Key words and phrases.} Chessboard complex, Determinantal ring,  Divisor class group, Finite free resolution, Koszul complex,  Segre ring,  Universal resolution.}

 \bigpagebreak
 
\flushpar{\bf   Introduction.}

\medskip

Let $G$ be a free module rank $g$ over the commutative noetherian ring $R$.
For each integer $n$, the complex 
$$\tsize \Bbb C_n:\quad \dots \to \operatorname{Sym}_{n-2}G\otimes{\tsize \bigwedge}^2G\to  \operatorname{Sym}_{n-1}G\otimes{\tsize \bigwedge}^1G\to \operatorname{Sym}_nG\otimes{\tsize \bigwedge}^0G\to 0$$ of free $R$-modules is well understood: the complex $\Bbb C_0$ has homology $\operatorname{H}_0(\Bbb C_0)=R$ and, for any  non-zero integer $n$, the complex $\Bbb C_n$ is split exact. 
 Indeed, $\Bbb C_n$  is one homogeneous  strand of a Koszul resolution. Let $x_1,\dots,x_g$ be a basis for $G$, then $\bigoplus_n\Bbb C_n$ is the minimal homogeneous resolution  of $R$ by free modules over the polynomial ring $R[x_1,\dots x_g]$. 

The situation is significantly different if two free modules $E^*$ and $G$ are used. In this case, the complexes look like 
 $$\eightpoint \dots \to \operatorname{Sym}_mE^*\otimes \operatorname{Sym}_nG\otimes {\tsize \bigwedge}^p(E^*\otimes G)\to \operatorname{Sym}_{m+1}E^*\otimes \operatorname{Sym}_{n+1}G\otimes {\tsize \bigwedge}^{p-1}(E^*\otimes G)\to \dots \tag{0.1}$$
and if $V\in \operatorname{Sym}_{m}E^*$, $X\in \operatorname{Sym}_{n}G$, $v_i\in E^*$ and $x_i\in G$, for $1\le i\le p$, then the differential 
$$\tsize \operatorname{Sym}_{m}E^*\otimes \operatorname{Sym}_{n}G\otimes {\tsize \bigwedge}^{p}(E^*\otimes G)\to  \operatorname{Sym}_{m+1}E^*\otimes \operatorname{Sym}_{n+1}G\otimes {\tsize \bigwedge}^{p-1}(E^*\otimes G)$$
sends $$V\otimes X\otimes (v_1\otimes x_1)\wedge (v_2\otimes x_2)\wedge \dots \wedge (v_{p}\otimes x_{p})$$to 
$$\sum_{i=1}^{p} (-1)^{i+1} v_iV\otimes x_iX\otimes (v_1\otimes x_1)\wedge  \dots \wedge \widehat{(v_i\otimes x_i)}\wedge \dots \wedge (v_{p}\otimes x_{p}).$$ Let $\operatorname{H}_{\Cal N}(m,n,p)$ equal the homology of ({0.1}) at $\operatorname{Sym}_{m}E^*\otimes \operatorname{Sym}_{n}G\otimes {\tsize \bigwedge}^{p}(E^*\otimes G)$. 

Many of the complexes ({0.1}) have homology. Sometimes the homology of ({0.1}) occurs someplace in the middle of ({0.1}) rather than only at one of the endpoints.  Some of the complexes ({0.1}) have homology in more than one position. The homology of ({0.1}) is not always free as an $R$-module. If $R$ is a field, then the dimension of the homology of ({0.1}) depends on the characteristic of $R$. 

We listed 5 differences between the complexes $\Bbb C_n$ and the complexes ({0.1}). The first three differences are not particularly surprising; however, the last two are stunning. These results were first  established by Hashimoto \cite{Ha}.  
Let $S$ be the $R$-algebra $\operatorname{Sym}_{\bullet}E^*\otimes \operatorname{Sym}_{\bullet}G$.
If 
 we fix bases $v_1,\dots,v_e$ for $E^*$, and $x_1,\dots,x_g$ for $G$, then one may think of $S$ as the polynomial ring 
$S=R[v_1, \dots, v_e, x_1, \dots, x_g]$. Let  $T$ be the subring $$T=\sum_m \operatorname{Sym}_mE^*\otimes \operatorname{Sym}_mG$$ of $S$. 
One may think of $T$ as the subring $R[v_ix_j]$ of $S$. Let $\Cal P$ be the R-algebra $\operatorname{Sym}_{\bullet}(E^*\otimes G)$. One may think of $\Cal P$ as a polynomial ring over $R$ in the $eg$ indeterminates  $\{v_i\otimes x_j\}$. It is convenient to let $z_{ij}$ represent the element $v_i\otimes x_j$ of $\Cal P$. The identity map on $E^*\otimes G$ induces a surjective map $\varphi\:\Cal P\to T$. Let $Z$ be the $e\times g$ matrix whose entry in row $i$ column $j$ is the indeterminate $z_{ij}$. The kernel of $\varphi$ is the ideal $I_2(Z)$ generated by the $2\times 2$ minors of $Z$; and therefore, $T$ is isomorphic to the determinantal ring $\Cal P/I_2(Z)$.  
Significant information about the homology of ({0.1}), in the case that $a=b$, is contained in   the graded  module $\operatorname{Tor}_{\bullet,\bullet}^{\Cal P}(T,R)$, where $R$ is the graded $\Cal P$-module $\Cal P/\Cal P_+$ concentrated in degree zero.    Hashimoto  proved that if $R$ is equal to the ring of integers, and $e$ and $g$ are both at least five, then $\operatorname{Tor}_{3,5}^{\Cal P}(T,R)$ is not a free $R$-module.
On the other hand, the Koszul complex $\Cal P\otimes_R {\tsize \bigwedge}^{\bullet}(E^*\otimes G)$ is a homogeneous resolution of the $\Cal P$-module $R$. It follows that 
$$\operatorname{Tor}^{\Cal P}_{p,n+p}(T,R)=\operatorname{H}_{\Cal N}(n,n,p).$$  The assertion about the dependence of the dimension of the homology of ({0.1}) on the characteristic of the base field follows immediately; see Roberts \cite{R}. 

There is a determinantal interpretation of the complexes ({0.1}), even when $m\neq n$. For each integer $\ell$, let $M_{\ell}$ be the $T$-submodule  
$$M_{\ell}=\sum\limits_{m-n=\ell} \operatorname{Sym}_mE^*\otimes \operatorname{Sym}_nG$$ of $S$. 
 View $M_{\ell}$ as a graded $T$-module by giving $\operatorname{Sym}_{n+\ell}E^*\otimes \operatorname{Sym}_n G$ grade $n$. The same reasoning we used before shows that  $$\operatorname{Tor}^{\Cal P}_{p,n+p}(M_{m-n},R)=\operatorname{H}_{\Cal N}(m,n,p).\tag{0.2}$$
The modules $M_{\ell}$ arise in numerous situations. For example, take $R=\Bbb Z$. The divisor class group of $T$ is known to be $\Bbb Z$ and \cite{BG} shows why $\ell\mapsto [M_\ell]$ is an isomorphism from $\Bbb Z\to \operatorname{C}\!\ell\, (T)$.  This numbering satisfies $M_0=T$, $M_{g-e}$ is equal to the canonical class of $T$, and $M_{\ell}$ is a Cohen-Macaulay $T$-module 
 if and only if $1-e\le \ell\le g-1$. 

Reiner and Roberts \cite{RR} refer to $T$ as the Segre ring $\operatorname{Segre(e,g,0)}$ and they call $M_{\ell}$ the Segre module $\operatorname{Segre}(e,g,\ell)$.

If $R=\pmb K$ is a field of characteristic zero, then the geometric method of Lascoux \cite{L} gives the modules in the resolution of $T$ by free $\Cal P$-modules and hence  dimension  of $\operatorname{Tor}^{\Cal P}_{p,q}(T,\pmb K)$ for all $p$ and $q$.  The same technique can be used to  resolve each $M_{\ell}$ over $\Cal P$; see \cite{W, Section 6.5}. A very pretty description of  
$\operatorname{Tor}^{\Cal P}_{p,q}(M_{\ell},\pmb K)$, when $\pmb K$ is a field of characteristic zero, is recorded in \cite{RR}; see section 1. The complete description of  
$\operatorname{Tor}^{\Cal P}_{p,q}(M_{\ell},\pmb K)$ for fields of arbitrary characteristic is not yet known.

The ring $\Cal P =\operatorname{Sym}_{\bullet}(E^*\otimes G)$ has a natural $\Bbb N^e\times \Bbb N^g$-grading; (once one picks bases for $E^*$ and $G$). Each module  $M_{\ell}$ and each homology module $\operatorname{Tor}^{\Cal P}_{p}(M_{\ell},\pmb K)$  inherits this grading.  It is shown in \cite{BH,RR}  that, for each $(\gamma,\delta)\in \Bbb N^e\times \Bbb N^g$, the $(\gamma,\delta)^{\text{th}}$ graded piece of  $\operatorname{Tor}^{\Cal P}_{p}(M_{\ell},\pmb K)$  is isomorphic to the $p-1$ reduced homology of the chessboard simplicial complex with multiplicities $\Delta_{\gamma,\delta}$. The vertex set of 
$\Delta_{\gamma,\delta}$ is the set of squares on an $e\times g$ chessboard; the simplices are the sets of squares having at most $\gamma_i$ squares from row $i$ and at most $\delta_j$ squares from column $j$ for all $i$ and $j$. The complexes  $\Delta_{\gamma,\delta}$ were introduced by Bruns and Herzog \cite{BH}; these complexes generalize the original chessboard complexes $\Delta_{e,g}$ of \cite{BLVZ}. The complex
$\Delta_{e,g}$ is equal to $\Delta_{\gamma,\delta}$ with $(\gamma,\delta)=((1,\dots,1),(1,\dots,1))$. 
It is shown in \cite{BLVZ} that $\widetilde{\operatorname{H}}_2(\Delta_{5,5},\Bbb Z)$ has $3$-torsion; thereby providing an independent proof of Hashimoto's Theorem. 
Bj\"orner,    Lov\'asz, Vre\'cica, and   \v Zivaljevi\'c make conjectures about the higher order connectivity of (or the higher order homotopy groups) of chessboard complexes. Friedman and Hanlon \cite{FH}  reformulate the conjectures of \cite{BLVZ} into statements about simplicial homology and they use characteristic zero techniques (Hodge Theory and eigenvalues of the Laplacian) to verify the conjectures in certain cases. 

The matching complex of a graph $\Gamma$ is the abstract simplicial complex whose vertex set is the set of edges 
of $\Gamma$ and whose faces are sets of edges of $\Gamma$ with no two edges meeting at a vertex. If $\Gamma$ is a complete bipartite graph, then the matching complex of $\Gamma$ is a chessboard complex. Karaguezian, Reiner, and Wachs \cite{KRW} use representation theory over fields of  characteristic zero to unite and extend all of the previously mentioned results. 

Almost all of the results listed above require a field of characteristic zero. We notice; however, that \cite{RR}, where the homology of ({0.1}) is explicitly recorded in characteristic zero, depends very heavily on duality. Duality continues to hold over $\Bbb Z$ -- see \cite{K05}.   
 Our goal in the present paper is to splice ({0.1}) together with the dual of the appropriate partner of ({0.1}) in order to produce a family of split exact complexes $\{\frak C^{r,s}\}$. This idea is motivated, for example, by the family of Eagon-Northcott like complexes (see, for example, \cite{E, Theorem A2.10})  or the family of complexes which resolves the divisor class group of a residual intersection of a grade three Gorenstein ideal \cite{KU}. Our approach, unlike that in \cite{RR, W, FH, KRW}, works over $\Bbb Z$. In other words, the same complex works in all characteristics. One side benefit of our approach is that we can pick out the honest-to-goodness cycles in ({0.1}) which generate the homology. 

The complex ({0.1})   and its dual play    prominent roles   in the resolution of a universal ring. For any triple of parameters $e$, $f$, and $g$, subject to the obvious constraints, Hochster \cite{Ho}  established the existence of a commutative noetherian ring $\Cal R$  and a universal resolution $$\Bbb U\:\quad 0\to \Cal R^{e}\to \Cal R^{f}\to \Cal R^{g}\to 0,\tag{0.3}$$ such that for any commutative noetherian ring $S$ and any resolution  $$\Bbb V\:\quad  0\to S^{e}\to S^{f}\to S^{g}\to 0,$$ there exists a unique ring homomorphism $\Cal R\to S$ with $\Bbb V=\Bbb U\otimes_{\Cal R} S$. In \cite{K06}, we found a free resolution $\Bbb F$ of  the universal ring $\Cal R$ over an integral polynomial ring, $\frak P$, in the border  case $f=e+g$. 
The resolution $\Bbb F$ is not minimal; indeed, if $e$ and $g$ are both at least $5$; then Hashimoto's Theorem  can be used to prove that $\Cal R$ does not posses a generic minimal resolution over the ring of integers. Nonetheless, for each field $\pmb K$, we are able to use the resolution $\Bbb F$ to express the the modules in the minimal homogeneous resolution of $\Cal R\otimes_{\Bbb Z}\pmb K$ by free $\frak P\otimes_{\Bbb Z} \pmb K$ modules in terms of the homology of ({0.1}). The homology of ({0.1}) is known when $\pmb K$ is a field of characteristic zero. With the answer in hand, but without appealing to 
\cite{K06}, the geometric method of calculating syzygies was directly applied  in \cite{KW} to find the modules in the minimal resolution of $\Cal R\otimes_{\Bbb Z}\pmb K$, when $\pmb K$ is a field of characteristic zero.

\medskip
 Let $E$ and $G$ be free modules of rank $e$ and $g$, respectively, over a commutative noetherian ring $R$. The identity map on $E^*\otimes G$ induces the Koszul complex ({0.1})  and its dual $$\eightpoint \dots \to D_{m+1}E\otimes D_{n+1}G^*\otimes {\tsize \bigwedge}^{p-1}(E\otimes G^*)\to D_mE\otimes D_nG^*\otimes {\tsize \bigwedge}^p(E\otimes G^*)\to \dots \ .\tag{0.4}$$Let  $\operatorname{H}_{\Cal M}(m,n,p)$ be  the homology of ({0.4}) at $D_mE\otimes D_nG^*\otimes {\tsize \bigwedge}^p(E\otimes G^*)$. 
  It is shown in \cite{BG}, implicitly, and \cite{K05}, explicitly, that 
$$\operatorname{H}_{\Cal N}(m,n,p)\cong \operatorname{H}_{\Cal M}(m',n',p'),\tag{0.5}$$
provided $m+m'=g-1$, $n+n'=e-1$, $p+p'={\pmb \alpha}$,  and
$1-e\le m-n\le g-1$, for ${\pmb \alpha}=(g-1)(e-1)$. 
 In Definition {2.3}  we define complexes $\Bbb N(P,Q)$ and $\Bbb M(P,Q)$, for integers $P$ and $Q$, in such a way that ({0.1}) is $\Bbb N(m+p,n+p)$, and ({0.4}) is $\Bbb M(m+p,n+p)$.   We exhibit  a quasi-isomorphism $$\CD \Bbb M(eg-e-P,eg-g-Q)[eg-P-Q]\\@V M VV\\ \Bbb N(P,Q),\endCD\tag{0.6}$$which yields an alternate proof of ({0.5}). The mapping cone of ({0.6}) is called ${\frak C}^{r,s}$ for $s=P$ and $r=e+P-Q$. 
We prove that ({0.6}) is a quasi-isomorphism by showing that each ${\frak C}^{r,s}$ is split exact.

In section 2 we define the objects  ${\frak C}^{r,s}$. We prove that each ${\frak C}^{r,s}$ is a complex in section 3; and that each  ${\frak C}^{r,s}$ is split exact in section  4. The most serious calculation in the paper is  the proof of Theorem {3.7}, where we produce a map $\psi\:{\frak C}^{r,s}\to{\frak C}^{r-1,s}$. We use $\psi$ to prove inductively that each ${\frak C}^{r,s}$ is a complex. We also use $\psi$ in section 4 to prove that each ${\frak C}^{r,s}$ is split exact. The map $\psi$ is very similar to the map $M$ of ({0.6}). The trick to the whole paper was to find these maps simultaneously. 
In section 5 we apply the quasi-isomorphism ({0.6}) in order  exhibit a generator of 
$$\operatorname{H}_{\Cal M}(g-1,e-1,{\pmb \alpha})\cong \operatorname{H}_{\Cal N}(0,0,0)=R.$$
The universal ring  $\Cal R$  of ({0.3}) is  a Gorenstein ring of projective dimension $eg+1$, and 
the module $\operatorname{H}_{\Cal M}(g-1,e-1,\pmb \alpha)$ contributes to the module in position $eg$ in the resolution of $\Cal R$.

In section 6 we use the technique of ({0.6}) to deduce a homogeneous version of ({0.5}). We also 
show that 
the homology  modules $\operatorname{H}_{\Cal M}(m,n,p)$ and $\operatorname{H}_{\Cal N}(m,n,p)$ satisfy an extra duality when $e$  is equal to three.


\bigpagebreak

\flushpar{\bf 1.\quad  The Reiner-Roberts description of $\operatorname{Tor}$ of a Segre module.}

\medskip

In Theorem {1.2} we record the homology $\operatorname{Tor}^{\Cal P}_{p,q}(M_{\ell},\pmb K)$,  when $\pmb K$ is a field of characteristic zero, as given in \cite{RR}. 
This description  is beautiful and is superior to other descriptions. 
In particular, if $\pmb K$ is a field of arbitrary characteristic, then the resolution of the universal ring $\Cal R\otimes_{\Bbb Z} \pmb K$ of ({0.3}) is given, in \cite{K06}, in terms of  the modules $\operatorname{Tor}^{\Cal P}_{p,q}(M_{\ell},\pmb K)$, with $-e\le \ell\le g-1$. On the other hand, if $\pmb K$ is a field of characteristic zero, then the resolution of $\Cal R\otimes_{\Bbb Z} \pmb K$ is given, in \cite{KW},  in terms of Schur modules. 
It is possible to combine the the results in \cite{K06} and \cite{KW} and describe  $\operatorname{Tor}^{\Cal P}_{p,q}(M_{\ell},\pmb K)$ in terms of Schur modules when $\pmb K$ has characteristic zero and $-e\le \ell\le g-1$. 
One can reparameterize to see that the description of $\operatorname{Tor}^{\Cal P}_{p,q}(M_{\ell},\pmb K)$ which is
 obtained from combining  \cite{K06} and \cite{KW} agrees with the description  given in \cite{RR}.

\definition{Data {1.1}} Let $E$ and $G$ be free modules of rank $e$ and $g$, respectively, over a commutative noetherian ring $R$, $\Cal P$ be the $R$-algebra $\operatorname{Sym}_{\bullet}(E^*\otimes G)$, $S$ be the $R$-algebra 
 $\operatorname{Sym}_{\bullet}E^*\otimes \operatorname{Sym}_{\bullet}G$, and   $T$ be the subring $$T=\sum_m \operatorname{Sym}_mE^*\otimes \operatorname{Sym}_mG$$ of $S$. 
For each integer $\ell$, let $M_{\ell}$ be the $T$-submodule  
$$M_{\ell}=\sum\limits_{m-n=\ell} \operatorname{Sym}_mE^*\otimes \operatorname{Sym}_nG$$ of $S$. 
 View $M_{\ell}$ as a graded $T$-module by giving $\operatorname{Sym}_{n+\ell}E^*\otimes \operatorname{Sym}_n G$ grade $n$. The identity map on $E^*\otimes G$ induces a surjective map $\varphi\:\Cal P\to T$; thereby making each $M_{\ell}$ a graded $\Cal P$-module. 
\enddefinition

Recall  the bi-graded structure on $\operatorname{Tor}$. 
\definition{Definition} If $\frak P=\bigoplus_i \frak P_i$ is a graded ring,   and $A=\bigoplus_i A_i$ and $B=\bigoplus_i B_i$ are graded $\frak P$-modules, then the module $\operatorname{Tor}^{\frak P}_{\bullet}(A,B)$ is a bi-graded $\frak P$-module. Indeed, if $$\Bbb Y\: \dots \to Y_1 \to Y_0\to A$$ is a $\frak P$-free  resolution of $A$,  homogeneous of degree zero, then 
$$\operatorname{Tor}_{p,q}^\frak P(A,B)=\frac{\ker [(Y_p\otimes B)_q\to (Y_{p-1}\otimes B)_q]}{\operatorname{im} [(Y_{p+1}\otimes B)_q\to (Y_{p}\otimes B)_q]}.$$\enddefinition

\proclaim{Theorem {1.2}}Adopt the notation of {1.1} with $R$ equal to a  field $\pmb K$ of characteristic zero.   For each integer $\ell$, 
$$\operatorname{Tor}^{\Cal P}_{\bullet,\bullet}(M_{\ell},{\pmb K})=\bigoplus\limits_{(\lambda,\mu)}S_{\lambda}E^*\otimes_{\pmb K} S_{\mu}G,$$
where $(\lambda,\mu)$ have the form 
$$\left(\matrix \format\l&\ \l\\ \boxed{{\phantom{XX}}}_{(s+1-\ell)\times s}&\beta\\\alpha\endmatrix\quad,\quad\matrix \format\l&\ \l\\ \boxed{{\phantom{XX}}}_{(s+1)\times (s-\ell)}&\alpha'\\  \beta'\endmatrix\right),$$ for some integer $s$ and partitions $\alpha$ and $\beta$ with
$s\ge \alpha_1$ and $\beta$ having at most $s-\ell$ parts.
Furthermore, the module $S_{\lambda}E^*\otimes_{\pmb K} S_{\mu}G$ is a summand of 
$$\operatorname{Tor}^{\Cal P}_{|\lambda|-s,|\lambda|-\ell}(M_{\ell},{\pmb K}).$$ The partitions $\lambda$ and $\mu$ have at most $e$ and $g$ parts respectively. 
\endproclaim

Our notation for the partition $\lambda$ represents what the Ferrers diagram for $\lambda$ looks like. In other words,  $\lambda_i=s+\beta_i$ for $1\le i\le s+1-\ell$ and $\lambda_{(s+1-\ell)+i}=\alpha_i$ for $1\le i$. We use the same notation for Schur modules as is used in \cite{KW}. 
The symbols $L_\lambda G$ and 
$K_\lambda G$ are used in \cite{W}  to denote the Schur and Weyl modules, respectively. In this paper we think about Schur modules only when we work 
over a field of characteristic zero; in which case,  our $S_\lambda G$
is isomorphic to $L_{\lambda^\prime}G$ or $K_\lambda G$, where
$\lambda^\prime$ is the conjugate partition.

The duality of ({0.5}) is the starting point of \cite{RR}; see especially Proposition 2.4 and Figure 3. When the base ring is a field $\pmb K$, we may use ({0.2}) to re-write ({0.5}) as
$$\operatorname{Tor}^{\Cal P}_{p,q}(M_{\ell},\pmb K) \cong \operatorname{Tor}^{\Cal P}_{p',q'}(M_{\ell'},\pmb K),\tag{1.3}$$
whenever $1-e\le \ell\le g-1$, $\ell+\ell'=g-e$,  $p+p'=(e-1)(g-1)$, $q+q'=(e-1)g$. It is somewhat remarkable that some notion of duality extends to the modules $\{M_{\ell}\}$ which are just outside the Cohen-Macaulay range. 
The complexes ${\frak C}^{r,s}$, with $r=0$ or $r=e+g$, of Definition {2.8} are the base case for our inductive arguments. These complexes were  introduced and shown to be split exact over any base ring in \cite{K05}.   
Some of the numerical consequences of the existence of these split exact complexes are captured in \cite{K05, Cor\. 2.2}.   Theorem {1.2} affords a quick proof of 
these numerical consequences (when the base ring is a field of characteristic zero).

\proclaim{Observation} If $\pmb K$ is a field of characteristic zero and the integers $p$ and $p'$ satisfy $p+p'=eg-e-g$, then 
$$\dim \operatorname{Tor}^{\Cal P}_{p,p+e}(M_{-e},\pmb K)+\dim \operatorname{Tor}^{\Cal P}_{p',p'}(M_{g},\pmb K)=\dim {\tsize \bigwedge}^{e+p}(E^*\otimes G).$$
\endproclaim
\demo{Proof} Theorem {1.2} gives
$$\matrix \operatorname{Tor}^{\Cal P}_{p,p+e}(M_{-e},\pmb K)=\sum\limits_{|\beta|=p+e\atop{\beta\subset e\times g\atop{\beta_1'=e}}}S_{\beta}E^*\otimes S_{\beta'}G\otimes {\tsize \bigwedge}^eE\ \text{and} \\\vspace{4pt} \operatorname{Tor}^{\Cal P}_{p,p}(M_g,\pmb K)=\sum\limits_{|\alpha|=p+g\atop{\alpha\subset e\times g\atop{\alpha_1=g}}}S_{\alpha}E^*\otimes S_{\alpha'}G\otimes {\tsize \bigwedge}^gG^*.\endmatrix
$$ 
We need two facts from Representation Theory; see, for example, \cite{W}:
$$S_{\lambda_1,\dots,\lambda_e}E^*=S_{-\lambda_e,\dots,-\lambda_1}E\quad \text{ and }\quad S_{\lambda_1,\dots,\lambda_e}E^*=S_{\lambda_1+1,\dots,\lambda_e+1}E^*\otimes {\tsize \bigwedge}^eE.$$
Fix $\alpha\subset e\times g$; let $\beta$ be the partition $g-\alpha_e,\dots,g-\alpha_1$. 
Notice that the Ferrers diagram of $\alpha$ fits together with the Ferrers diagram of $\beta$ (after it has been rotated by $180$ degrees) to form an $e\times g$ rectangle. In particular,  $\beta\subset e\times g$; $|\alpha|+|\beta|=eg$;   $\alpha_1=g\iff \beta_1'<e$; and $\beta'$ is the partition $e-\alpha_g',\dots,e-\alpha_1'$. We see that 
$$S_{\alpha}E^*=S_{\beta}E\otimes({\tsize \bigwedge}^eE^*)^{\otimes g}\quad\text{and}\quad S_{\alpha'}G=S_{\beta'}G^*\otimes({\tsize \bigwedge}^gG)^{\otimes e}.$$ 
If $p+p'=eg-e-g$, then  $$\operatorname{Tor}^{\Cal P}_{p',p'}(M_g,\pmb K)=\sum\limits_{|\alpha|=p'+g\atop{\alpha\subset e\times g\atop{\alpha_1=g}}}S_{\alpha}E^*\otimes S_{\alpha'}G\otimes {\tsize \bigwedge}^gG^*$$
$$=\sum\limits_{|\beta|=p+e\atop{\beta\subset e\times g\atop{\beta_1'<e}}}S_{\beta}E\otimes S_{\beta'}G^*\otimes {\tsize \bigwedge}^gG^*\otimes({\tsize \bigwedge}^gG)^{\otimes e}\otimes({\tsize \bigwedge}^eE^*)^{\otimes g}.$$
It follows that 
$$\dim \operatorname{Tor}^{\Cal P}_{p,p+e}(M_{-e},\pmb K)+\dim \operatorname{Tor}^{\Cal P}_{p',p'}(M_g,\pmb K)$$$$
=\sum\limits_{|\beta|=p+e\atop{\beta\subset e\times g\atop{\beta_1'=e}}}\dim S_{\beta}E^*\otimes S_{\beta'}G
+\sum\limits_{|\beta|=p+e\atop{\beta\subset e\times g\atop{\beta_1'<e}}}\dim S_{\beta}E\otimes S_{\beta'}G^*$$
$$=\sum\limits_{|\beta|=p+e\atop{\beta\subset e\times g}}\dim S_{\beta}E^*\otimes S_{\beta'}G=\dim {\tsize \bigwedge}^{p+e} (E^*\otimes G).$$The final equality is the Cauchy formula. \qed 
\enddemo


\bigpagebreak

\flushpar{\bf 2.\quad  The complexes ${\frak C}^{r,s}$.}

\medskip

Throughout this paper, $R$ is a commutative noetherian ring with one, $E$ and $G$ are fixed, non-zero, finitely generated  free $R$-modules of rank $e$ and $g$, respectively, and $${\pmb \alpha}=(e-1)(g-1).$$  Fix bases $$\text{ $x_1,\dots,x_g$ for $G$; $y_1,\dots,y_g$ for $G^*$; $u_1,\dots,u_e$ for $E$; and $v_1,\dots,v_e$ for $E^*$,}\tag{2.1}$$ with $\{x_i\}$ dual to $\{y_j\}$, and $\{u_i\}$ dual to $\{v_j\}$. 

\definition{Definition {2.2}} Let   $m$, $n$, and $p$ be integers. Define modules
 $$\split  \Cal M(m,n,p)={}& D_mE\otimes D_nG^*\otimes {\tsize {\tsize \bigwedge}}^p(E\otimes G^*)\ \text{and}\\ 
   \Cal N(m,n,p)={}& \operatorname{Sym}_mE^*\otimes \operatorname{Sym}_{n}G\otimes {\tsize {\tsize \bigwedge}}^{p}(E^*\otimes G),\endsplit$$
and homomorphisms
  $$\split   \pmb D\:\Cal M(m,n,p)&{}\to \Cal M(m-1,n-1,p+1)\ \text{and}\\
 K\: \Cal N(m,n,p)&{} \to \Cal N(m+1,n+1,p-1).\\
\endsplit
$$In the language of ({2.1}), 
$$\split \pmb D(U\otimes   Y\otimes Z){}&=\sum v_k(U)\otimes x_{\ell}(Y)\otimes (u_k\otimes y_{\ell})\wedge Z,\ \text{and}\\
 K(V\otimes   X\otimes W){}&=\sum v_kV\otimes x_{\ell}X\otimes (u_k\otimes y_{\ell})(W).\endsplit$$
Each sum is taken over all $k$ and $\ell$, with $1\le k\le e$ and $1\le \ell \le g$.
 \enddefinition

 \definition{Definition {2.3}}  Fix integers $P$ and $Q$. Let $\Bbb N(P,Q)$ 
be the complex 
$$0\to \Cal N(P-eg,Q-eg,eg)@>K>> \dots @>K>> \Cal N(P-1,Q-1,1)@>K>> \Cal N(P,Q,0)\to 0,$$
and $\Bbb M(P,Q)$ be 
$$0\to \Cal M(P,Q,0)@>\pmb D>> \Cal M(P-1,Q-1,1)@>\pmb D>> \dots@>\pmb D>>  \Cal M(P-eg,Q-eg,eg)\to 0.$$
The module $\Cal N(P,Q,0)$ is in position zero in $\Bbb N(P,Q)$; $\Cal M(P,Q,0)$ is in position $P+Q+1$ in $\Bbb M(P,Q)$.  \enddefinition

The  perfect paring  $$\Cal M(m,n,p)\otimes \Cal N(m,n,p)\to R,\tag{2.4}$$ which is given by   
$$(U\otimes   Y\otimes Z)\otimes (V\otimes X\otimes    W)\mapsto U(V)\cdot Y(X)\cdot Z(W),$$ensures that the maps  $\pmb D$ and $K$ are dual to one another. The details are given below.

\proclaim{Lemma {2.5}}
If $T\in \Cal M(m,n,p)$ and $S\in \Cal N(m-1,n-1,p+1)$, then 
 $$[K(S)](T)=(-1)^{p}S[\pmb D(T)].$$
\endproclaim
\demo{Proof}Let $T=U\otimes Y\otimes Z$ and $S=V\otimes X\otimes W$. We use the notation of {(\rm}{2.1}{\rm)}. The left side is equal to 
$$\split &\sum [v_k\cdot V\otimes x_{\ell}\cdot X\otimes (u_k\otimes y_{\ell})(W)][U\otimes Y\otimes Z]\\&{}=\sum [v_k\cdot V](U)\cdot [x_{\ell}\cdot X](Y)\cdot [(u_k\otimes y_{\ell})(W)](Z)\\&{}=(-1)^p\sum V(v_k(U))\cdot   X (x_{\ell}(Y))\cdot W[(u_k\otimes y_{\ell})\wedge Z],\endsplit$$ and this is the right side of the proposed identity. Each sum is taken over all $k$ and $\ell$, with $1\le k\le e$ and $1\le \ell \le g$.  \qed
\enddemo

\definition{Notation}{\bf (a)} Let $m$ be an integer. Each pair of elements $(U,Y)$, with $U\in D_{m}E$ and  $Y\in {\tsize \bigwedge}^{m}G^*$,  gives rise to an element of ${\tsize \bigwedge}^m(E\otimes G^*)$, which we denote by 
$U\bowtie Y$. We now give the definition of $U\bowtie Y$. Consider the composition 
$$ D_mE\otimes  T_mG^*@>\Delta\otimes 1>> T_mE\otimes  T_mG^*  @>\xi >>   {\tsize \bigwedge}^m(E\otimes G^*), $$where 
$ \xi\left((U_1\otimes \dots \otimes U_m)\otimes (Y_1\otimes\dots\otimes Y_m) {\vphantom{E^{E^{E}}_{E_{E}}}} \right)= (U_1\otimes Y_1)\wedge \dots \wedge (U_m\otimes Y_m)$,
for $U_i\in E$ and $Y_i\in G^*$.  It is easy to see  that the above composition factors through $D_mE\otimes {\tsize \bigwedge}^mG^*$. Let $U\otimes Y\mapsto U\bowtie Y$ be the resulting map from $D_mE\otimes {\tsize \bigwedge}^mG^*$ to ${\tsize \bigwedge}^m(E\otimes G^*)$. The map 
$${\tsize \bigwedge}^mE\otimes D_mG^*\to {\tsize \bigwedge}^m(E\otimes G^*),$$which sends $U\otimes Y$ to $U\bowtie Y$, for $U\in {\tsize \bigwedge}^mE$ and $Y\in  D_mG^*$, is defined in a completely analogous manner.
\flushpar  {\bf (b)} For each  statement  ``S'', let $$\chi(\text{S})= \cases 1,&\text{if S is true, and}\\ 0,&\text{if S is false.}\endcases$$In particular, $\chi(i=j)$ has the same value as the Kronecker delta $\delta_{ij}$.

\flushpar  {\bf(c)} Every free $R$-module that we consider is oriented; in the sense that,  if $F$ is a free module of rank $f$, then $\omega_F$ is the name of our preferred generator for ${\tsize \bigwedge}^fF$. The orientations of $F$ and $F^*$ are always compatible in the sense that $\omega_F(\omega_{F^*})=1$ and $\omega_{F^*}(\omega_F)=1$.
\enddefinition

In Definition {2.7} we introduce the complexes ${\frak C}^{r,s}$. The modules and the part of the differential which is independent of the choice of coordinates is found in Definition {2.2}. The part of the differential which  is coordinate dependent is  called ``$M_m$'' and is introduced in Definition {2.6}. The maps ``$N_m$'' are very similar to the $M_m$. The maps ``$N_m$'' are not part of ${\frak C}^{r,s}$, but are used in Theorem {3.7} to produce a map of complexes ${\frak C}^{r,s}\to {\frak C}^{r-1,s}$. 

\definition{Definition}Adopt the notation of {(\rm}{2.1}{\rm)}. Fix  an integer $j$,   with $1\le j\le g$. 
\flushpar{\bf (a)} Let $Y_j^{-}$ and $Y_j^{+}$ be the elements $$Y_j^{-}=y_1\wedge\dots \wedge y_j\in {\tsize \bigwedge}^jG^*\quad\text{and}\quad   Y_j^{+}=y_j\wedge\dots \wedge y_g\in {\tsize \bigwedge}^{g-j+1}G^*.$$Let $Y_0^-=Y^+_{g+1}=1$. 
\flushpar{\bf (b)} Define     the $R$-module homomorphism $\tau_j\:D_{\bullet}G^* \otimes D_{\bullet}G^* \to D_{\bullet}G^*$    by   
$$\tau_j(y_1^{(a_1)}\cdots y_g^{(a_g)}\otimes y_1^{(b_1)}\cdots y_g^{(b_g)})=\chi\cdot y_1^{(a_1)}\cdots y_{j-1}^{(a_{j-1})}y_j^{(a_j+b_j+1)}y_{j+1}^{(b_{j+1})}\cdots y_g^{(b_g)},  $$ where $\chi=\chi(0\le a_j)\chi(0\le b_j)\chi(b_1=\dots= b_{j-1}=a_{j+1}=\dots=a_g=0)$.
\enddefinition

\definition{Definition {2.6}}Adopt the notation of {(\rm}{2.1}{\rm)}. Let $\omega=\omega_{E^*\otimes G}$. For integers $m$, $n$, and $p$, we   define   maps $$\matrix \format \l\\
M_m\:\Cal M(m,n,p)\to\Cal N(g-1-m,e-1-n,{\pmb \alpha}-p) \ \text{and}\\
N_m\:\Cal M(m,n,p)\to\Cal N(g-2-m,e-1-n,{\pmb \alpha}+1-p).
\endmatrix$$Fix   $T=U\otimes Y\otimes Z\in \Cal M(m,n,p)$.  If $T'=U'\otimes Y'\otimes Z'\in\Cal M(g-1-m,e-1-n,{\pmb \alpha}-p)$, then 
 $[M_m(T)](T')$ is equal to 
$$\left[ \left(\omega_E\bowtie \tau_{m+1}(Y'\otimes Y)\right)\wedge  Z'\wedge  Z \wedge \left(U\bowtie Y^{-}_{m} \right)
 \wedge\left(U' \bowtie Y^{+}_{m+2}\right) \right](\omega).
$$
If $T'=U'\otimes Y'\otimes Z'\in\Cal M(g-2-m,e-1-n,{\pmb \alpha}+1-p)$, then  
$[N_m(T)](T')$  is equal to $$ (-1)^{m+p}\left[ \left(\omega_E\bowtie \tau_{m+1}( Y' \otimes Y)\right) \wedge   Z'\wedge  Z \wedge \left(U\bowtie  Y^{-}_{m} \right)
\wedge\left(U'\bowtie (x_g(Y^{+}_{m+2}))\right) \right](\omega).$$
 \enddefinition

   \definition{Definition {2.7}}For integers $r$ and $s$, with  $1\le r\le e+g-1$, we define  the complex $${\frak C}^{r,s}\:\quad 
  \cdots \to {\frak C}_i^{r,s}@>d_i>> {\frak C}_{i-1}^{r,s}\to \cdots ,$$
by$${\frak C}_i^{r,s}= \Cal M(g+i-s-2,r+i-s-2,{\pmb \alpha}-i+1) \oplus   \Cal N(s-i,e+s-r-i,i)$$ and $$d_i=\left[\matrix\format\l&\ \ \l\\ \pmb D&0\\M_{g+i-s-2}&K\endmatrix\right].$$
\enddefinition
 
\remark{Remark}
If $1\le r\le e+g-1$, then ${\frak C}^{r,s}$ is the mapping cone of the following map of complexes:
$$\eightpoint \CD
 \dots  @>\pmb D>>\Cal M(g-s,r-s,{\pmb \alpha}-1)@>\pmb D >> \Cal M(g-s-1,r-s-1,{\pmb \alpha})@>\pmb D >>\dots\\  {}@. @V -M_{g-s} VV @V M_{g-s-1} VV\\ \dots   @> K>>  \Cal N(s-1,e+s-r-1,1) @> K>> \Cal N(s,e+s-r,0)@>>> 0,\endCD$$with $\Cal N(s,e+s-r,0)$ is position zero. This description of ${\frak C}^{r,s}$ was promised in ({0.6}).
\endremark

Complexes ${\frak C}^{r,s}$, with $r=0$ or $r=e+g$ 
 have already been introduced and shown to be split exact in \cite{K05}. These complexes form the base case for our inductive arguments and  we recall them at this point. 
\definition{Definition {2.8}}{\bf (a)} For each integer $p$, let $B(p)={\tsize \bigwedge}^p(E\otimes G^*)$  and define 
 homomorphisms
  $$\split   \gamma\:\Cal M(0,e,p)&{} \to B(e+p),\\
\gamma\:\Cal M(g,0,p)&{} \to B(g+p),\\
\Gamma\: B(p)&{} \to \Cal N(g,0,eg-g-p),\ \text{and}\\
\Gamma\: B(p)&{} \to \Cal N(0,e,eg-e-p).
\endsplit
$$The maps $\gamma$ are described by $$\gamma(1\otimes Y\otimes Z)=   (\omega_E \bowtie Y)\wedge Z  \quad\text{and}\quad  \gamma(U\otimes 1\otimes Z)=   (U \bowtie \omega_{G^*})\wedge Z.$$
If $Z$ is in $B(p)$, then $\Gamma(Z)$ is the element of $\Cal N(g,0,eg-g-p)$  (or  $\Cal N(0,e,eg-e-p)$, respectively) with
$$[\Gamma(Z)](T)=[\gamma(T)\wedge Z](\omega_{E^*\otimes G})$$ for all 
$T$ in $\Cal M(g,0,eg-g-p)$ (or $\Cal M(0,e,eg-e-p)$, respectively). 
 
\flushpar{\bf (b)} If $r=0$, then ${\frak C}^{0,s}$ is 
$$\split 0\to \Cal M(eg-e-s,{\pmb \alpha}-s-1,0)@>\pmb D>> \dots@>\pmb D>>\Cal M(g,0,{\pmb \alpha}-s-1)@>\gamma>> B(eg-e-s)\\ @>\Gamma>> \Cal N(0,e,s)@>K>> \dots @>K>> \Cal N(s,e+s,0)\to 0,\endsplit $$ with ${\frak C}_0^{0,s}=\Cal N(s,e+s,0)$.

\flushpar{\bf (c)} If $r=e+g$, then ${\frak C}^{e+g,s}$ is 
$$\nopagebreak\split 0\to \Cal M(eg-e-s,eg-s,0)@>\pmb D>> \dots@>\pmb D>>\Cal M(0,e,eg-e-s)@>\gamma>> B(eg-s)\\ @>\Gamma>> \Cal N(g,0,s-g)@>K>> \dots @>K>> \Cal N(s,s-g,0)\to 0,\endsplit $$ with ${\frak C}_0^{e+g,s}=\Cal N(s,s-g,0)$.
\enddefinition

\remark{Remarks {2.9}}
{\bf (a)} The description of ${\frak C}^{r,s}_i$ which is given in Definition {2.7} holds for $r=0$, provided $i\neq s+1$. This description also holds for $r=e+g$, provided $i\neq s-g+1$. Remark {3.10}  explains why it is necessary to include the $B$ summands in ${\frak C}^{0,s}$ and ${\frak C}^{e+g,s}$. 

\flushpar{\bf (b)} Fix $r$, with   $0\le r\le e+g$. 
 Observe that $\Cal N(m,n,p)$ is a summand of ${\frak C}^{r,s}_i$ if $$r=e+m-n,\quad s=m+p,\quad\text{and}\quad  i=p;$$ and
$\Cal M(m,n,p)$ is a summand of ${\frak C}^{r,s}_i$ if $$r=g-m+n,\quad s=eg-e-m-p,\quad\text{and}\quad  i={\pmb \alpha}-p+1.$$

\flushpar{\bf (c)} If $$r+r'=e+g, \quad s+s'=eg-e,\quad\text{and}\quad   i+i'={\pmb \alpha}+1,$$ then ({2.4}) provides  a natural perfect pairing 
${\frak C}^{r,s}_i\otimes  {\frak C}^{r',s'}_{i'}\to R$.
\endremark

\example{Example {2.10}}Take $g=1$ and $0\le r\le e+g$. In this case, the complex ${\frak C}^{r,s}$ is a Koszul complex or the  dual  of a Koszul complex. In fact,  ${\frak C}^{r,s}$ is equal to  
$$\alignat 2& 0\to \operatorname{Sym}_0E^*\otimes{\tsize \bigwedge}^sE^*@>K>>\dots@>K>>\operatorname{Sym}_sE^*\otimes {\tsize \bigwedge}^0E^*\to 0, &&\qquad\text{if $1\le s$},\\
 &0\to R@>= >> R\to 0,&&\qquad \text{if $s=0$},\\
 &0\to D_{-s}E\otimes {\tsize \bigwedge}^{0}E@>\pmb D >> \dots @>\pmb D >> D_{0}E\otimes {\tsize \bigwedge}^{-s}E \to 0,  &&\qquad\text{if $s\le -1$}. \endalignat  $$  All of these complexes are   split exact. 
In the language of Definitions {2.7} and {2.8}, the complex ${\frak C}^{r,s}$ is
$$\alignat 2
&0\to \Cal N(0,e-r,s)\to \dots\to \Cal N(s,e-r+s,0)\to 0,&&\ \text{if $1\le s$,  $r\le e$},\\
&0\to B(e-s)\to \Cal N(1,0,s-1)\to \dots\to \Cal N(s,s-1,0)\to 0,&&\ \text{if $1\le s$,  $r=e+1$},\\
&0\to B(0)\to \Cal N(0,e,0)\to 0,&&\ \text{if $s=0$,  $r=0$},\\
&0\to \Cal M(0,r-1,0)\to \Cal N(0,e-r,0)\to 0,&&\ \text{if $s=0$,  $1\le r\le e$},\\
&0\to \Cal M(0,e,0)\to B(e)\to 0,&&\ \text{if $s=0$,  $r=e+1$},\\
&0\to \Cal M(-s,-s-1,0)\to\hskip-.14pt \dots \hskip-.14pt \to \Cal M(1,0,-s-1)\to B(-s)\to 0,&&\ \text{if $s\le -1$,  $r=0$},\\
&0\to \Cal M(-s,r-1-s,0)\to \dots \to \Cal M(0,r-1,-s)\to 0,&&\ \text{if $s\le -1$,  $1\le r$}.\endalignat$$
The modules $\Cal N(s,e-r+s,0)$ and $\Cal M(-s,r-1-s,0)$ are in positions zero and one, respectively.  \endexample

\example{Example {2.11}} The complex ${\frak C}^{1,{\pmb \alpha}}$ is
$$0\to \Cal M(g-1,0,0)@>M_{g-1}>> \Cal N(0,e-1,{\pmb \alpha})@>K>> \dots @>K>> \Cal N({\pmb \alpha},eg-g,0)\to 0.$$The map $M_{g-1}$ depends on the choice of coordinates far less than appears to be the case. For each $j$, let $\int\underline{\phantom{X}}\, d\!y_j\: D_{\bullet}G^*\to D_{\bullet}G^*$ be the homomorphism which sends $y_1^{(a_1)}\cdots y_g^{(a_g)}$ to $y_1^{(a_1)}\cdots y_j^{(a_j+1)}\cdots y_g^{(a_g)}$. Observe that $$[M_{g-1}(U\otimes 1\otimes 1)][1\otimes Y'\otimes Z']=\pm \left[\left(\omega_{E}\bowtie  \tsize \int Y' \,d\!y_g\right)\wedge \left(U\bowtie x_g(\omega_{G^*})\right)\wedge Z'\right](\omega).$$
Corollary {3.3}~(b), applied to $\int\!\!\int Y'\, d\!y_j\, d\!y_g$,  shows that the elements
$$\tsize (\omega_{E}\bowtie \int Y' \,d\!y_j)\wedge (U\bowtie x_j(\omega_{G^*})) $$
of $\bigwedge^{e+g-1}(E\otimes G^*)$ are equal for all $j$, with $1\le j\le g$. 
\endexample

\example{Example}  The complex ${\frak C}^{2,{\pmb \alpha}+1}$ is
$$\split {0\to {\matrix \Cal M(g-2,0,0)\\\oplus\\ \Cal N(0,e-2,{\pmb \alpha}+1)\endmatrix } @>{\bmatrix M_{g-2}&K\endbmatrix }>> \Cal N(1,e-1,{\pmb \alpha})@>K>> \dots} &\\ {\dots @>K>> \Cal N({\pmb \alpha}+1,eg-g,0)\to 0.}&\endsplit $$
Let $U$ be an element of $D_{g-2}E$, $\hat{G}$ be the free submodule of $G$ generated by $x_1,\dots,x_{g-1}$, and $\hat{M}_m$ be the map of Definition {2.6} constructed with data $(E,\hat{G})$. Example {2.11} shows that $$\hat{M}_{g-2}(U\otimes 1\otimes 1)$$ is an unbounded cycle in 
${\Cal N(0,e-1,(e-1)(g-2))}$.  It is clear that 
$$\xi=(1\otimes y_g\otimes 1)\left(K(1\otimes 1\otimes (\omega_{E^*}\bowtie x_g^{(e)}))\right)$$is an unbounded cycle in $\Cal N(1,0,e-1)$.  It is not difficult to see that $M_{g-2}(U\otimes 1\otimes 1)$ is equal to the product $\pm \hat{M}_{g-2}(U\otimes 1\otimes 1)\times \xi$ in $\Cal N(1,e-1,{\pmb \alpha})$.
\endexample

\proclaim{Observation} If $r$, $r'$, $s$, and $s'$ are integers with $0\le r\le e+g$, $r+r'=e+g$, and $s+s'=eg-e$, then the complexes ${\frak C}^{r,s}$ and $({\frak C}^{r',s'})^*[-({\pmb \alpha}+1)]$ are isomorphic. \endproclaim
\demo{Proof}Remark {2.9}~(c) shows that the corresponding modules are isomorphic. Lem\-ma {2.5} records the duality between $K$ and $\pmb D$. The maps $\gamma$ and $\Gamma$ are defined to be dual to one another. 
Construct the maps ``$M^{\text{rev}}_m$'' of Definition {2.6} using the reversed order $y_g,\dots,y_1$ for the basis of $G^*$. It is clear that $[M_m^{\text{rev}}(T) ](T')$ is equal to $\pm [M_{g-1-m}(T')](T)$  for $T\in \Cal M(m,n,p)$ and $T'\in \Cal M (g-1-m,e-1-n,{\pmb \alpha}-p)$.
 \qed\enddemo

The next result is well-known. One proof of it   appears in \cite{K06}. The result shows that 
one may prove a formula about the elements of $D_{m}E$  by checking that the formula holds at each pure divided power $u^{(m)}$, for $u\in E$, provided  the formula can be obtained, by way of base change, from a corresponding formula over a  polynomial ring over the ring of integers. 
 \proclaim{Lemma {2.12}} Suppose $R$ is a polynomial ring over the ring of integers,   $E$ and $G$ are free $R$-modules, and $\varphi\:D_mE\to G$ is an $R$-module homomorphism. If $\varphi(u^{(m)})=0$ for all $u\in E$, then $\varphi$ is identically zero.\endproclaim

 \proclaim{Observation {2.13}}Adopt the notation of {\rm({2.1})}. If $m$ is an integer, $y\in G^*$, $Y$ is in ${\tsize \bigwedge}^mG^*$, and $U\in D_{m+1}E$, then
$$U\bowtie (y\wedge Y)=\sum\limits_{k=1}^e (u_k\otimes y)\wedge (v_k(U)\bowtie Y).$$\endproclaim
\demo{Proof}Each side represents a homomorphism from $D_{m+1}E$ to ${\tsize \bigwedge}^{m+1}(E\otimes G^*)$. It suffices to check that equality holds for $U=u^{(m+1)}$ with  $u\in E$. \qed\enddemo


\bigpagebreak

 \flushpar{\bf 3.\quad   Each ${\frak C}^{r,s}$ is a complex.}

\medskip

In this section we prove Theorem {3.1}.
 \proclaim{Theorem {3.1}} Let $r$ and $s$ be integers. 
\flushpar{\bf(a)} If $0\le r\le e+g$, then   ${\frak C}^{r,s}$ is a complex.  
\flushpar{\bf(b)} If $1\le r\le e+g-1$, then 
$\psi\:{\frak C}^{r,s}\to {\frak C}^{r-1,s}$,  from Definition {3.4},  is a map of complexes. 
 \endproclaim

\demo{Proof}Corollary {3.3}~(a) contains the only interesting step in the proof  that  each ${\frak C}^{0,s}$ is a complex. 
If $1\le r\le e+g-1$, then we exhibit a commutative diagram ${\frak C}^{r,s}\to {\frak C}^{r-1,s}$ in Observation {3.5} and Theorem {3.7}. It is clear that  $K$ and $\pmb D$  square to zero. We need only worry about the part of $d\circ d$ which involve the maps $M$ of Definition {2.6}. The maps $M$ land in the $\bigoplus \Cal N(m,n,p)$ part of ${\frak C}^{r,s}$ and the map ${\frak C}^{r,s}\to {\frak C}^{r-1,s}$ is obviously injective on the $\Cal N(m,n,p)$ summands of ${\frak C}^{r,s}$. Induction on $r$ yields that each ${\frak C}^{r,s}$ is a complex. \qed \enddemo

Lemma {3.2}, which appears to be rather innocuous, is actually quite powerful. We have recorded a few of the elementary consequences of this result in Corollary {3.3}.  We have already appealed to Corollary {3.3}~(b) in Example {2.11}. Also, part (a) of  {3.3} appears in \cite{K05}, where it is used to prove that  ${\frak C}^{r,s}$ is  a complex for $r$ equal to $0$ or $e+g$.

 \proclaim{Lemma {3.2}}Adopt the notation of {\rm({2.1})}. 
If $m$ is an integer,  ${Y\in D_{e+1}G^*}$, $U$ is in $D_{m+1}E$,  and ${y\in {\tsize \bigwedge}^{m}G^*}$, then
$$\sum_{\ell=1}^{g}  \left(\omega_E\bowtie [x_{\ell}(Y)]\right)\wedge \left(U\bowtie [y_{\ell}\wedge y]\right)$$is equal to zero in ${\tsize \bigwedge}^{e+m+1}(E\otimes G^*)$. 
\endproclaim
 \demo{Proof}In light of Lemma {2.12}, we may assume that the base ring is a domain and that $U=u^{(m+1)}$ for some non-zero $u\in E$. There exists $v\in E^*$, with $u(v)$ not equal to  zero. Use the fact that 
$$u\wedge v(\omega_E) =[u(v)](\omega_E) $$for $u\in {\tsize \bigwedge}^{\bullet}E$ and $v\in {\tsize \bigwedge}^{\bullet}E^*$, as well as Observation {2.13}, to see that  $u(v)$  times the indicated  expression is 
$$\sum_{\ell} \left([u\wedge v(\omega_E)]\bowtie [x_{\ell}(Y)]\right)\wedge \left(u^{(m+1)}\bowtie [y_{\ell}\wedge y]\right)$$
$$=\sum_{k,\ell} (u\otimes y_k)\wedge\left(v(\omega_E)\bowtie [(x_kx_{\ell})(Y)]\right)\wedge \left(u^{(m+1)}\bowtie [y_{\ell}\wedge y]\right)$$
$$=(-1)^{e-1} \sum_{k,\ell}   \left(v(\omega_E)\bowtie [(x_kx_{\ell})(Y)]\right)\wedge \left(u^{(m+2)}\bowtie [y_k\wedge y_{\ell}\wedge y]\right),$$ and this  is zero because $\sum\limits_{k,\ell}x_kx_{\ell}\otimes y_k\wedge y_{\ell}$ is zero in $\operatorname{Sym}_2G \otimes {\tsize \bigwedge}^2G^*$. \qed
\enddemo

\proclaim{Corollary {3.3}}  
{\bf (a)} If    $Y\in D_{e}G^*$ and $U\in D_{g}E$,  then
$$ \nopagebreak (\omega_E\bowtie Y)\wedge (U\bowtie \omega_{G^*})$$is equal to zero in ${\tsize \bigwedge}^{e+g}(E\otimes G^*)$. 
\flushpar{\bf (b)} If    $Y\in D_{e+1}G^*$, $x$ and $x'$ are in $G$, and $U\in D_{g-1}E$,  then
$$ \nopagebreak (\omega_E\bowtie x(Y))\wedge (U\bowtie x'(\omega_{G^*}))\quad\text{and}\quad  (\omega_E\bowtie x'(Y))\wedge (U\bowtie x(\omega_{G^*}))$$  are equal in ${\tsize \bigwedge}^{e+g-1}(E\otimes G^*)$. 
\flushpar{\bf (c)} In the notation of {\rm({2.1})}, if    $Y\in D_{e}G^*$, $U\in D_{g-1}E$, and $x_j(Y)=0$, for some $j$, with $1\le j\le g$,     then
$$ \nopagebreak (\omega_E\bowtie  Y )\wedge (U\bowtie x_j(\omega_{G^*}))=0$$    ${\tsize \bigwedge}^{e+g-1}(E\otimes G^*)$. 
\endproclaim
 \demo{Proof}Take $Y'\in D_{e+1}G^*$ with $y_1(Y')=Y$. Observe that 
$$(\omega_E\bowtie Y)\wedge (U\bowtie \omega_{G^*})=\sum_{\ell} (\omega_E\bowtie x_{\ell}(Y'))
\wedge(U\bowtie [y_{\ell}\wedge x_1(\omega_{G^*})])
,$$and this is zero by Lemma {3.2}.  
For (b), apply Lemma {3.2} with $y=(x\wedge x')(\omega_{G^*})$. 
 In (c), take $j$ to be $g$, for   notational convenience. The assertion is  true, but not interesting, unless $2\le g$. We may assume that $Y$ is a basis vector, say $y_1^{(a_1)}\cdots y_g^{(a_g)}$. Let $Y'$ be a basis vector with $x_1(Y')=Y$ and $x_g(Y')=0$. Apply (b).  
\qed
\enddemo

\definition{Definition {3.4}}Let $r$ and $s$ be integers.
If $r=1$, then define $\psi\:{\frak C}^{1,s}\to {\frak C}^{0,s}$ by
$$\hskip-3pt\eightpoint \CD \cdots@>\pmb D>>\Cal M(g,1,{\pmb \alpha}-s-1)@>\pmb D>> \Cal M(g-1,0,{\pmb \alpha}-s)  @>M_{g-1}>> \Cal N(0,e-1,s)@>K>> \cdots   \\
@. @V \pm x_g VV @V \nu VV @V x_g VV @.\\
\cdots@>\pmb D>>\Cal M(g,0,{\pmb \alpha}-s-1)@>\gamma>> B(eg-e-s)  @>\Gamma>> \Cal N(0,e,s)@>K>> \cdots. 
\endCD$$
If $2\le r\le e+g-1$, then define 
$\psi_i\:{\frak C}^{r,s}_i\to {\frak C}^{r-1,s}_i$ to be $$\bmatrix \format\l&\ \ \l\\ x_g&0\\N_{g+i-s-2}&x_g\endbmatrix.$$ The  map $x_g\:\Cal M(m,n,p)\to \Cal M(m,n-1,p)$ is given by the action of $x_g$ on $D_nG^*$;  
$x_g\:\Cal N(m,n,p)\to \Cal N(m,n+1,p)$ is given by multiplication by $x_g$ in $\operatorname{Sym}_{\bullet}G$; and if $T=U\otimes 1\otimes Z\in \Cal M(g-1,0,{\pmb \alpha}-s)$, then $\nu(T)=Z\wedge (U\bowtie Y_{g-1}^{-})$.   The module ${\frak C}^{r,s}_i$ is given in Definition {2.7} and the map $N$ is given in Definition {2.6}.\enddefinition

\proclaim{Observation {3.5}}For each integer $s$, $\psi\: {\frak C}^{1,s}\to {\frak C}^{0,s}$  is a commutative diagram. \endproclaim
\demo{Proof}Consider the diagram of Definition {3.4}. Let $T=U\otimes 1\otimes Z\in \Cal M(g-1,0,{\pmb \alpha}-s)$ and $T'=1\otimes Y'\otimes Z'\in \Cal M(0,e,s)$. We see that $$[(\Gamma\circ \nu)(T)](T')=[(\omega_E\bowtie Y')\wedge Z'\wedge Z\wedge (U\bowtie Y_{g-1}^{-})](\omega),$$ and $$[(1\otimes x_g\otimes 1)\times M_{g-1}(T)](T')=\left[\left(\omega_E\bowtie [\tau_g(x_g(Y')\otimes 1)]\right)\wedge Z'\wedge Z\wedge (U\bowtie Y_{g-1}^{-})\right](\omega).$$ It is clear that $x_g$ sends $\tau_g(x_g(Y')\otimes 1)-Y'$ to zero. Use Corollary {3.3}c to conclude that $\Gamma\circ \nu= (1\otimes x_g\otimes 1)\times M_{g-1}$. Now let $T=u^{(g)}\otimes y\otimes Z\in \Cal M(g,1,{\pmb \alpha}-s-1)$. We see that $$[\nu\circ \pmb D](T)= 
(u\otimes y)\wedge Z\wedge (u^{(g-1)}\bowtie Y_{g-1}^{-})$$ and $$[\gamma\circ  (1\otimes x_g\otimes 1)](T)=  x_g(y)\cdot (u^{(g)}\bowtie \omega_{G^*})\wedge Z;$$thus, $\pm \gamma\circ  (1\otimes x_g\otimes 1)= \nu\circ \pmb D$, for the appropriate choice of sign. 
\qed \enddemo

The homomorphisms $L_j$ and $\Upsilon_j$ are projections. For example, $L_j$ projects onto the summand of $D_{\bullet}G^*$ which does not involve  any basis vector with a subscript less than or equal to $j$. 
\definition{Definition}Adopt the notation of {(\rm}{2.1}{\rm)}. Fix  an integer $j$,   with $1\le j\le g$. 
 Define $R$-module homomorphisms  $L_{j}\: D_{\bullet}G^* \to D_{\bullet}G^*$  and $\Upsilon_j\: D_{\bullet}G^* \to D_{\bullet}G^*$   by $$\matrix \format\l\\L_j(y_1^{(b_1)}\cdots y_g^{(b_g)})=\chi(b_1=\dots=b_j=0)\cdot y_1^{(b_1)}\cdots y_g^{(b_g)}\text{ and}\\ \Upsilon_j(y_1^{(b_1)}\cdots y_g^{(b_g)})=\chi(b_j=\dots=b_g=0)\cdot y_1^{(b_1)}\cdots y_g^{(b_g)}.\endmatrix $$Define  $L_0=\Upsilon_{g+1}$   to be  the identity map. 
\enddefinition

\proclaim{Lemma {3.6}}Fix    integers $j$ and $\ell$,   with $1\le j,\ell\le g$. Let $Y$ and $Y'$ be elements of $D_{\bullet}G^*$. The following statements hold:
 \flushpar{\bf  (a)}
 if $\ell<j$, then $\tau_j(x_{\ell}(Y')\otimes Y)=x_{\ell}(\tau_j(Y'\otimes Y))$,
\flushpar{\bf  (b)} $x_{j}(\tau_j(Y'\otimes Y))=\tau_j(x_{j}(Y')\otimes Y)+\Upsilon_j(Y')\cdot L_{j-1}(Y)$,
\flushpar{\bf  (c)}
 if $j<\ell$, then $\tau_j(Y'\otimes x_{\ell}(Y))=x_{\ell}(\tau_j(Y'\otimes Y))$, and 
\flushpar{\bf  (d)} $x_{j}(\tau_j(Y'\otimes Y))=\tau_j(Y'\otimes x_{j}(Y))+\Upsilon_{j+1}(Y')\cdot L_{j}(Y)$. \endproclaim
\demo{Proof}It suffices to assume $Y'=y_1^{(a_1)}\cdots y_g^{(a_g)}$ and $Y=y_1^{(b_1)}\cdots y_g^{(b_g)}$. We prove (b). Let $Y''$ be the element
$$Y''=\chi\cdot  y_1^{(a_1)}\cdots y_{j-1}^{(a_{j-1})}y_j^{(a_j+b_j)}y_{j+1}^{(b_{j+1})}\cdots y_g^{(b_g)}$$
of $D_{\bullet}G^*$, where $\chi=\chi(0\le b_j)\chi(b_1=\dots= b_{j-1}=a_{j+1}=\dots=a_g=0)$.
Observe  that 
$$\matrix \format\l&\l\\ 
x_{j}(\tau_j(Y'\otimes Y))&{}=\chi(0\le a_j)Y'',\\\vspace{5pt}
\tau_j(x_{j}(Y')\otimes Y)&{}=\chi(1\le a_j)Y'', \ \text{and}\\\vspace{5pt}
\Upsilon_j(Y')\cdot L_{j-1}(Y)&{}=\chi(0= a_j)Y''.
\endmatrix$$The bottom formula  holds, even for $j=1$, since  $L_{0}$ is the identity map. The proof of (b) is complete. The proof of (d) is similar. Assertions (a) and (c) are obvious. \qed \enddemo

\proclaim{Theorem {3.7}}If $r$ and $s$ are integers, with $2\le r\le e+g-1$, then $$\psi\:{\frak C}^{r,s}\to {\frak C}^{r-1,s}$$ is a commutative diagram.
 \endproclaim

\demo{Proof}
Let $m$, $n$, and $p$ be integers and with $r=g-m+n$. We will prove that  
 the square
$$\nopagebreak \CD 
\Cal M(m,n,p)@>{\left[\smallmatrix \pmb D\\M_m\endsmallmatrix\right]}>>
{\matrix   \Cal M(m-1,n-1,p+1)\\ \oplus\\  \Cal N(g-1-m,e-1-n,{\pmb \alpha}-p)\endmatrix}\\
@V{\left[\smallmatrix x_g\\N_m\endsmallmatrix\right]}VV @V{\left[\smallmatrix N_{m-1}&x_g \endsmallmatrix\right]}VV\\
{\matrix \Cal M(m,n-1,p)\\\oplus\\ \Cal N(g-2-m,e-1-n,{\pmb \alpha}+1-p)\endmatrix}@>{\left[\smallmatrix M_{m}&K \endsmallmatrix\right]}>>  \Cal N(g-1-m,e-n,{\pmb \alpha}-p)
\endCD
$$commutes. We may as well assume that $0\le m\le g-1$ and 
$0\le n\le e$; otherwise, there is nothing to prove. 
Let  $T=U\otimes Y\otimes Z$ in $\Cal M(m,n,p)$,  $T'=U'\otimes Y'\otimes Z'$ in ${\Cal M(g-1-m,e-n,{\pmb \alpha}-p)}$, $z=Z'\wedge Z\in {\tsize \bigwedge}^{{\pmb \alpha}}(E\otimes G^*)$,  $$ \matrix \format\l&\quad \l\\
\Cal A=\chi(1\le m) [(N_{m-1}\circ \pmb D)(T)](T'),&
\Cal B=[(1\otimes x_g\otimes 1)\times M_{m}(T)](T'), \\
\Cal C=[M_{m}((1\otimes x_g\otimes 1)(T))](T'),\ \text{and}&   
\Cal D=[(K\circ N_{m})(T)](T'). \endmatrix$$
We will show that $\Cal A+\Cal B=\Cal C+\Cal D$.   Observe that $\Cal A$ is 
$$ \ssize  \chi(1\le m) \sum\limits_{k,\ell}
(-1)^{m}\left[ \left(\omega_E\bowtie \tau_{m}( Y' \otimes x_{\ell}(Y))\right) \wedge   z \wedge  (u_k\otimes y_{\ell})\wedge \left(v_k(U)\bowtie  Y^{-}_{m-1} \right)
\wedge\left(U'\bowtie (x_g(Y^{+}_{m+1}))\right) \right](\omega).$$We know that $\chi(1\le m)\neq 0$ implies $0<m$; hence, Observation {2.13} shows that $\Cal A$ is equal to $$ \ssize    \chi(1\le m)\sum\limits_{\ell}
(-1)^{m}\left[ \left(\omega_E\bowtie \tau_{m}( Y' \otimes x_{\ell}(Y))\right) \wedge   z \wedge \left(U\bowtie  [y_{\ell}\wedge Y^{-}_{m-1}] \right)
\wedge\left(U'\bowtie (x_g(Y^{+}_{m+1}))\right) \right](\omega).$$The expression $y_{\ell}\wedge Y^{-}_{m-1}$ is zero unless $m\le \ell\le g$; so,  we may apply Lemma {3.6} to see that $\Cal A=\Cal A_1+\Cal A_2$, with $\Cal A_1$ equal to 
$$\ssize  \chi(1\le m)\sum\limits_{m\le \ell\le g}
(-1)^{m}\left[ \left(\omega_E\bowtie x_{\ell}(\tau_{m}( Y' \otimes Y))\right) \wedge   z \wedge \left(U\bowtie  [y_{\ell}\wedge Y^{-}_{m-1}] \right)
\wedge\left(U'\bowtie (x_g(Y^{+}_{m+1}))\right) \right](\omega)$$
and $ \Cal A_2 $ equal to $$\ssize \chi(1\le m)
(-1)^{m+1}\left[ \left(\omega_E\bowtie \Upsilon_{m+1}(Y')L_{m}(Y)\right) \wedge   z \wedge \left(U\bowtie  [y_{m}\wedge Y^{-}_{m-1}] \right)
\wedge\left(U'\bowtie (x_g(Y^{+}_{m+1}))\right) \right](\omega).$$The value of $\Cal A_1$ is not changed if the index $\ell$ runs from $1$ to $g$; thus, Lemma {3.2} tells us that $\Cal A_1=0$.  Notice that in $\Cal A_2$,
$$(-1)^{m+1}y_{m}\wedge Y^{-}_{m-1}= Y^{-}_{m}.$$ 
The map $\Upsilon_1$ is the identity map on $D_0G^*$ and the zero map on $D_qG^*$ for all positive $q$.
It follows that  the factor $\chi(1\le m)$, in $\Cal A_2$, is redundant, because this factor only removes the case $m=0$. However, if $m$ were to equal $0$, then  $\Upsilon_{m+1}(Y')$ would already be 
$$\Upsilon_1(Y')=\chi(n=e)Y'=\chi(r=e+g)Y',\tag{3.8}$$ and this is zero, because the ambient hypothesis requires $r\le e+g-1$. We conclude  that  
$$  
\Cal A=\left[ \left(\omega_E\bowtie \Upsilon_{m+1}(Y')L_{m}(Y)\right) \wedge   z \wedge \left(U\bowtie Y^{-}_{m} \right)
\wedge\left(U'\bowtie (x_g(Y^{+}_{m+1}))\right) \right](\omega).$$
We see that $$  \Cal B=\left[ \left(\omega_E\bowtie \tau_{m+1}(x_g(Y')\otimes Y)\right)\wedge  z \wedge \left(U\bowtie Y^{-}_{m} \right)
 \wedge\left(U' \bowtie Y^{+}_{m+2}\right) \right](\omega)$$and
$$ \Cal C =\left[ \left(\omega_E\bowtie \tau_{m+1}(Y'\otimes x_g(Y))\right)\wedge  z \wedge \left(U\bowtie Y^{-}_{m} \right)
 \wedge\left(U' \bowtie Y^{+}_{m+2}\right) \right](\omega).
$$
Use Lemma {2.5} to see that 
$ \Cal D =(-1)^{{\pmb \alpha}-p}[N_{m}(T)][\pmb D(T')]$ is equal to  
 $$  \ssize    \sum\limits_{k,\ell} (-1)^{m+{\pmb \alpha}}\left[ \left(\omega_E\bowtie \tau_{m+1}( x_{\ell}(Y') \otimes Y)\right) \wedge   (u_k\otimes y_{\ell}) \wedge z \wedge \left(U\bowtie  Y^{-}_{m} \right)
\wedge\left(v_k(U')\bowtie (x_g(Y^{+}_{m+2}))\right) \right](\omega).$$If $m=g-1$, then $U'\in D_0E$ and $\Cal D$ is zero. Otherwise, we apply  Observation {2.13} to see that    $\Cal D$ is equal to 
 $$  \eightpoint   \chi(m\le g-2)\sum\limits_{\ell}\left[ \left(\omega_E\bowtie \tau_{m+1}( x_{\ell}(Y') \otimes Y)\right) \wedge    z \wedge \left(U\bowtie  Y^{-}_{m} \right)
\wedge\left(U'\bowtie (y_{\ell}\wedge x_g(Y^{+}_{m+2}))\right) \right](\omega).$$
The expression $y_{\ell}\wedge x_g(Y^{+}_{m+2})$ is zero, unless $\ell\le m+1$ or $\ell=g$. We write $\Cal D=\Cal D_1+\Cal D_2$, with $ \Cal D_1$ equal to  
$$  \chi(m\le g-2)\left[ \left(\omega_E\bowtie \tau_{m+1}( x_{g}(Y') \otimes Y)\right) \wedge    z \wedge \left(U\bowtie  Y^{-}_{m} \right)
\wedge\left(U'\bowtie Y^{+}_{m+2}\right) \right](\omega)$$ and $\Cal D_2$ equal to 
$$   \ssize  \chi(m\le g-2)\sum\limits_{1\le \ell\le m+1}\left[ \left(\omega_E\bowtie \tau_{m+1}( x_{\ell}(Y') \otimes Y)\right) \wedge    z \wedge \left(U\bowtie  Y^{-}_{m} \right)
\wedge\left(U'\bowtie (y_{\ell}\wedge x_g(Y^{+}_{m+2}))\right) \right](\omega).$$
Apply Lemma {3.6} to write $\Cal D_2=\Cal D_2'+\Cal D_2''$, where $\Cal D_2'$ is
$$ \ssize    \chi(m\le g-2)\sum\limits_{1\le \ell\le m+1}\left[ \left(\omega_E\bowtie x_{\ell}(\tau_{m+1}( Y' \otimes Y))\right) \wedge    z \wedge \left(U\bowtie  Y^{-}_{m} \right)
\wedge\left(U'\bowtie (y_{\ell}\wedge x_g(Y^{+}_{m+2}))\right) \right](\omega)$$
and $ \Cal D_2''$ is 
$$    \chi(m\le g-2)\left[ \left(\omega_E\bowtie \Upsilon_{m+1}(Y')L_{m}(Y)\right) \wedge    z \wedge \left(U\bowtie  Y^{-}_{m} \right)
\wedge\left(U'\bowtie (x_g(Y^{+}_{m+1}))\right) \right](\omega).$$
Apply Lemmas {3.2} and  {3.6} to see that  $\Cal D_2'$ is equal to  $$    -\chi(m\le g-2)\left[ \left(\omega_E\bowtie \tau_{m+1}( Y' \otimes x_{g}(Y))\right) \wedge    z \wedge \left(U\bowtie  Y^{-}_{m} \right)
\wedge\left(U'\bowtie Y^{+}_{m+2}\right) \right](\omega).$$
Observe that  $\Cal A-\Cal D_2''$ is equal to 
 $$   \chi(m=g-1)\left[ \left(\omega_E\bowtie \Upsilon_{g}(Y')L_{g-1}(Y)\right) \wedge    z \wedge \left(U\bowtie  Y^{-}_{g-1} \right)
\wedge\left(U'\bowtie 1\right) \right](\omega),$$
and that $\Cal B-\Cal D_1$ is equal to 
$$ \chi(m=g-1)\left[ \left(\omega_E\bowtie \tau_{g}( x_{g}(Y') \otimes Y)\right) \wedge    z \wedge \left(U\bowtie  Y^{-}_{g-1} \right)
\wedge\left(U'\bowtie 1\right) \right](\omega).$$Use Lemma {3.6} to see that $\Cal A+\Cal B-(\Cal D_1+\Cal D_2'')$ is equal to 
$$   \chi(m=g-1)\left[ \left(\omega_E\bowtie x_{g}(\tau_{g}( Y' \otimes Y))\right) \wedge    z \wedge \left(U\bowtie  Y^{-}_{g-1} \right)
\wedge\left(U'\bowtie 1\right) \right](\omega).$$
Use Lemma {3.6}, again, to express $x_{g}(\tau_{g}( Y' \otimes Y))$ as $\tau_{g}( Y' \otimes x_{g}(Y))+\Upsilon_{g+1}(Y')L_g(Y)$. The map  $L_g$ acts like  the identity   on $D_0G^*$ and like zero  on $D_qG^*$ for all positive $q$. Thus,  $$\chi(m=g-1)L_g(Y)=\chi(m=g-1)\chi(n=0)Y,\tag{3.9}$$ and this product is zero because the ambient hypothesis guarantees that $r$, which is equal to $g-m+n$, is at least $2$. Thus, $\Cal A+\Cal B-(\Cal D_1+\Cal D_2'')$ is equal to $$   \chi(m=g-1)\left[ \left(\omega_E\bowtie \tau_{g}( Y' \otimes x_{g}(Y))\right) \wedge    z \wedge \left(U\bowtie  Y^{-}_{g-1} \right)
\wedge\left(U'\bowtie 1\right) \right](\omega),$$which is equal to  $\Cal C+\Cal D_2'$, and the proof is complete. \qed
\enddemo
\remark{Remark {3.10}} Our proof of   Theorem {3.7} does not work for  either $r=1$  or $r=e+g$. Indeed, neither of  the  squares     
$$\CD 
\Cal M(g-1,0,p)@> M_{g-1} >>
  \Cal N(0,e-1,{\pmb \alpha}-p) \\
@V VV @V x_g VV\\
 0@> >>  \Cal N(0,e,{\pmb \alpha}-p),
\endCD
$$ nor 
$$\CD 
\Cal M(0,e,p)@>   >>
  0\\
@V x_g VV @VVV\\
 \Cal M(0,e-1,p)@> M_0  >>  \Cal N(g-1,0,{\pmb \alpha}-p)
\endCD
$$commute.  The proof breaks down for $m=0$ and $r=e+g$ at ({3.8}). It breaks down for $m=g-1$ and $r=1$ at ({3.9}). It is fortunate that our complex 
${\frak C}^{0,{\pmb \alpha}-p}$ is 
$$\cdots \to \Cal M(g,0,p-1)\to B(g+p-1)\to \Cal N(0,e,{\pmb \alpha}-p)\to \cdots $$ rather than 
  $$\cdots \to \Cal M(g,0,p-1)\to \Cal M(g-1,-1,p)\to\Cal N(0,e,{\pmb \alpha}-p)\to \cdots;$$ and that ${\frak C}^{e+g,eg-e-p}$ is 
$$\cdots \to \Cal M(0,e,p)\to B(e+p)\to \Cal N(g,0,{\pmb \alpha}-p-1)\to \cdots$$ rather than $$\cdots\to \Cal M(0,e,p)\to
 \Cal N(g-1,-1,{\pmb \alpha}-p)\to \Cal N(g,0,{\pmb \alpha}-p-1)\to \cdots.$$
\endremark


\bigpagebreak

\flushpar{\bf 4.\quad  Exactness.}

\medskip

\proclaim{Theorem {4.1}} If $r$ and $s$ are integers with $0\le r\le e+g$, then the complex ${\frak C}^{r,s}$ is split exact.
\endproclaim
\demo{Proof}It is shown in \cite{K05, Cor\. 2.26} show that the result holds for $r=0$ or $r=e+g$. Example {2.10} takes care of the case $g=1$. The proof proceeds by induction on $r$ and $g$. Fix integers $r$ and $s$, with $1\le r\le e+g-1$. Let  $\psi\:{\frak C}^{r,s}\to {\frak C}^{r-1,s}$ be the map of complexes  from Theorem {3.1}, and let   $(\Bbb A,\pmb a)$ be the mapping cone of $\psi$. We know, by induction, that ${\frak C}^{r-1,s}$ is split exact. It will suffice to show that $\Bbb A$ is exact, and this isn't very difficult. We split off two split exact subcomplexes from $\Bbb A$ in order to produce the complex ``$\Bbb B/\Bbb P$'', which has the same homology as $\Bbb A$. Then we show that $\Bbb B/\Bbb P$ is isomorphic to the direct sum of complexes which is given in ({4.2}). The complexes of ({4.2}) are known to be split exact because they are made using $g-1$ in place of $g$ and $r-1$ in place of $r$.

 The complex  $\Bbb A$ looks like 
$$\dots\to \Bbb A_i@> \pmb a_i>> \Bbb A_{i-1} \to \dots, $$ with $\Bbb A_i= {\frak C}^{r,s}_{i-1}\oplus {\frak C}^{r-1,s}_i$ and 
$$\pmb a_i=\bmatrix d_{i-1}&0\\ \psi_{i-1}&-d_i \endbmatrix.$$  Recall the bases of ({2.1}). Separate $x_1,\dots,x_{g-1}$ from $x_g$, and $y_1,\dots,y_{g-1}$ from $y_g$. Use $D_nG^*=D_n^0G^*\oplus D_n^+G^*$ and $\operatorname{Sym}_nG=\operatorname{Sym}_n^0G\oplus \operatorname{Sym}_n^+G$ to obtain $$\split \Cal M(m,n,p)={}&\Cal M^0(m,n,p)\oplus \Cal M^+(m,n,p)\ \text{and}\\ \Cal N(m,n,p)={}&\Cal N^0(m,n,p)\oplus \Cal N^+(m,n,p),\endsplit $$ 
 where the basis elements of $D_n^0G^*$ and $\operatorname{Sym}_n^0G$ involve $y_g^{(0)}$ and $x_g^0$, respectively; and the basis elements of $D_n^+G^*$ and $\operatorname{Sym}_n^+G$ involve $y_g^{(p)}$ and $x_g^{p}$, respectively, for positive powers $p$.
Let $\Bbb B$ be the submodule of $\Bbb A$  which consists of all modules of the form 
$\Cal M^0(m,n,p)$ or $\Cal N(m,n,p)$ from ${\frak C}^{r,s}$ and all modules of the form $B(m)$  or $\Cal N(m,n,p)$ from ${\frak C}^{r-1,s}$. It is easy to see that $\pmb a(\Bbb B)\subseteq \Bbb B$; hence, $(\Bbb B,\pmb a)$ is a subcomplex of $(\Bbb A,\pmb a)$. The quotient $\Bbb A/\Bbb B$ consists of the modules 
$\Cal M^+(m,n,p)$  from ${\frak C}^{r,s}$ and  $\Cal M(m,n,p)$ from ${\frak C}^{r-1,s}$, with differential 
$$ \dots@>>>
{\matrix   \Cal M^+(m,n,p)\\\oplus\\  \Cal M(m+1,n,p-1)
\endmatrix}
@>{\left[\matrix
\operatorname{proj}\circ \pmb D&0\\
x_g&-\pmb D
\endmatrix \right]}
>> 
{\matrix 
 \Cal M^+(m-1, n-1, p+1)\\\oplus\\  \Cal M(m,n-1,p)\endmatrix}@>>> \dots\ . $$ 
The map $x_g$ is an isomorphism;
 hence, $\Bbb A/\Bbb B$ is exact and $\operatorname{H}_i(\Bbb A)=\operatorname{H}_i(\Bbb B)$ for all $i$.

Let $\Bbb P$ be the submodule of $\Bbb B$  which consists of all modules of the form 
$\Cal N(m,n,p)$ from ${\frak C}^{r,s}$ and all modules of the form $\Cal N^+(m,n,p)$ from ${\frak C}^{r-1,s}$. It is easy to see that $(\Bbb P,\pmb a)$ is a subcomplex of $(\Bbb B,\pmb a)$; indeed, 
 $\Bbb P$ looks like 
$$ \dots@>>>
{\matrix \Cal N(m,n,p)
\\\oplus\\ \Cal N^+(m-1,n,p+1)
\endmatrix}
@>{\left[\matrix
K&0\\
x_g&-K
\endmatrix \right]}
>> 
{\matrix 
 \Cal N(m+1,n+1,p-1)\\\oplus\\
 \Cal N^+(m,n+1,p)\endmatrix}@>>> \dots \ .$$ 
 Each map $x_g$ is an isomorphism; thus, $\Bbb P$ is split exact and   $\operatorname{H}_i(\Bbb A)=\operatorname{H}_i(\Bbb B/\Bbb P)$ for all $i$. 
If $r=1$, then 
the complex  $ \Bbb B/\Bbb P $ is 
$$\eightpoint  \cdots@>\pmb D>>\Cal M^0(g,1,{\pmb \alpha}-s-1)@>\pmb D>> \Cal M^0(g-1,0,{\pmb \alpha}-s)  @> \nu  >> B(eg-e-s)  @>\operatorname{proj}\circ \Gamma>> \Cal N^0(0,e,s)@>\operatorname{proj}\circ K>> \cdots. 
$$
If $2\le r$, then the module  $(\Bbb B/\Bbb P)_i$ is 
$$ \matrix  \Cal M^0(m,n,p)\\\oplus\\ \Cal N^0(g-3-m,e-2-n,{\pmb \alpha}+2-p),
\endmatrix$$with $\Cal M^0(m,n,p)$ a summand of ${\frak C}_{i-1}^{r,s}$ and $\Cal N^0(g-3-m,e-2-n,{\pmb \alpha}+2-p)$
a summand of ${\frak C}_{i}^{r-1,s}$. (One may use Definition {2.7} to express $m$, $n$, and $p$, in terms of $r$, $s$, and $i$.) The differential $ (\Bbb B/\Bbb P)_i\to (\Bbb B/\Bbb P)_{i-1}$ is 
$${\bmatrix \pmb D&0\\\operatorname{proj}\circ N_m&-\operatorname{proj}\circ K\endbmatrix}.$$

Let $\hat{G}$ and ${\hat{G}}^*$ be the free submodules of $G$ and $G^*$ which are generated by $x_1,\dots,x_{g-1}$ and $y_1,\dots, y_{g-1}$, respectively. Let $\hat{\Cal M}(m,n,p)$, $\hat{\Cal N}(m,n,p)$, $\hat{M}_m$ and $({\hat{{\frak C}}}^{r,s},\hat{d})$ be the modules, maps,  and complexes manufactured using the data $E$ and $\hat{G}$. We will  show that the complex $\Bbb B/\Bbb P$ is is isomorphic to the following direct sum of complexes:
$$\bigoplus_{\ell=0}^e({\hat{{\frak C}}}^{r-1,s+\ell-e}[\ell-e]\otimes {\tsize {\tsize \bigwedge}}^{\ell}E, \hat{d}\otimes 1).\tag{4.2}$$The proof will then be complete, since each of the above complexes   is exact by induction on $r$. 

Consider the homomorphisms 
$$\split \Phi\: {}&\tsize \bigoplus_{\ell}\hat{\Cal M}(m,n,p-\ell)\otimes {\tsize {\tsize \bigwedge}}^{\ell}E\to \Cal M^0(m,n,p)\ \text{and}\\
\Phi\: {}&\tsize \bigoplus_{\ell}\hat{\Cal N}(m,n,p-\ell)\otimes {\tsize {\tsize \bigwedge}}^{\ell}E^*\to \Cal N^0(m,n,p),\endsplit $$which are given by $$\split\Phi((U\otimes Y\otimes Z)\otimes u)={}&U\otimes Y\otimes Z\wedge (u\bowtie y_g^{(\ell)})\ \text{and}\\
\Phi((V\otimes X\otimes W)\otimes v)={}&V\otimes X\otimes W\wedge (v\bowtie x_g^{(\ell)}),\endsplit$$for
$U\otimes Y\otimes Z\in \hat{\Cal M}(m,n,p-\ell)$, $u\in  {\tsize {\tsize \bigwedge}}^{\ell}E$, $V\otimes X\otimes W\in \hat{\Cal N}(m,n,p-\ell)$, and $v\in {\tsize {\tsize \bigwedge}}^{\ell}E^*$. It is clear that each map $\Phi$ is an isomorphism; and therefore,   $\Phi$ induces a module isomorphism $\pmb \Phi$ from ({4.2}) to $\Bbb B/\Bbb P$. Lemma {4.3} takes care of the only tricky part of showing that $\pmb \Phi$ is an isomorphism of complexes.
 \qed \enddemo

\proclaim{Lemma {4.3}} Retain the notation from the proof of Theorem {\rm{4.1}}.
Fix integers $k$ and $\ell$. Consider the composition
$$\split \hat{\Cal M}(m,n,p-\ell)\otimes {}&{}{\tsize \bigwedge}^{\ell}E@>\Phi>> \Cal M^0(m,n,p)\\{}&{}@>\operatorname{proj}\circ N_m>> \Cal N^0(g-2-m,e-1-n,{\pmb \alpha}+1-p)\\{}&{}
@>\operatorname{proj}\circ \Phi^{-1}>> \hat{\Cal N}(g-2-m,e-1-n,{\pmb \alpha}+1-p-k)\otimes {\tsize \bigwedge}^{k}E^*\\{}&{}
@>>> \hat{\Cal N}(g-2-m,e-1-n,{\pmb \alpha}+1-p-k)\otimes {\tsize \bigwedge}^{e-k}E,\endsplit $$ where the last map sends $v\in {\tsize \bigwedge}^kE^*$ to $v(\omega_E)$. 
If $\ell+k=e$, then the above composition is $\pm \hat{M}_m\otimes 1$; otherwise, the composition is zero. 
\endproclaim
\demo{Proof}There is nothing to prove unless $m\le g-2$. Let $T=U\otimes Y\otimes Z$ in $\hat{\Cal M}(m,n,p-\ell)$ and $u\in {\tsize \bigwedge}^{\ell}E$. We compare the elements 
$$\nopagebreak \operatorname{proj}\circ N_m\circ \Phi(T\otimes u)\quad\text{and}\quad  \chi(\ell+k=e)\Phi(\hat{M}_m(T)\otimes u(\omega_{E^*}))$$ of the submodule $$\nopagebreak  \Phi(\hat{\Cal N}(g-2-m,e-1-n,{\pmb \alpha}+1-p-k)\otimes {\tsize \bigwedge}^{k}E^*) $$ of $\Cal N^0(g-2-m,e-1-n,{\pmb \alpha}+1-p)$. 
In other words, we compare 
$$\split \Cal A={}& [N_m\circ \Phi(T\otimes u)][\Phi(T'\otimes u')]\ \text{and}\\  \Cal B={}& \chi(\ell+k=e)[\Phi(\hat{M}_m(T)\otimes u(\omega_{E^*}))][\Phi(T'\otimes u')],\endsplit $$
for $T'=U'\otimes Y'\otimes Z'\in \hat{\Cal M}(g-2-m,e-1-n,{\pmb \alpha}+1-p-k)$ and $ u' \in {\tsize \bigwedge}^kE$. 
Use Definition {2.6} to expand 
$$\Cal A= [N_m(U\otimes Y\otimes Z\wedge (u\bowtie y_g^{(\ell)}))][U'\otimes Y'\otimes Z'\wedge (u'\bowtie y_g^{(k)})].$$
In the expanded version of $\Cal A$, all of  the $y_g$'s  appear  in the expressions $(u\bowtie y_g^{(\ell)})$ and $(u'\bowtie y_g^{(k)})$ and  all of the $x_g$'s appear as the factor $\omega_{E^*}\bowtie x_g^{(e)}$ of $\omega$. Let $\hat{\omega}$ represent $\omega_{E^*\otimes \hat{G}}$. It follows that $\Cal A$ is equal to 
$$C\cdot  \left[\left(\omega_E\bowtie \tau_{m+1}(Y'\otimes Y)\right)\wedge Z'\wedge  Z\wedge \left(U\bowtie Y_m^-\right) \wedge \left(U'\bowtie (x_g(Y^+_{m+2}))\right)\right](\hat{\omega}),$$ where $C$ is the constant $ \pm \chi(\ell+k=e)(u'\wedge u)(\omega_{E^*})$. One quickly sees that $\Cal B$ equals $C\cdot [\hat{M}_m(T)][T']$.
Furthermore, $[\hat{M}_m(T)][T']$ is equal to 
$$ \left[ \left(\omega_E\bowtie \hat{\tau}_{m+1}(Y'\otimes Y)\right)\wedge  Z'\wedge  Z \wedge \left(U\bowtie {\hat{Y}}^{-}_{m} \right)
 \wedge\left(U' \bowtie {\hat{Y}}^{+}_{m+2}\right) \right](\hat{\omega}),
$$ where the map $\hat{\tau}_{m+1}$ and the elements ${\hat{Y}}^{-}_{m}$ and  ${\hat{Y}}^{+}_{m+2}$ are manufactured using the data $E$ and $\hat{G}$. Neither $Y$ nor $Y'$ involve $y_g$ and $m+1\le g-1$; so, $\hat{\tau}_{m+1}(Y'\otimes Y)$ equals $\tau_{m+1}(Y'\otimes Y)$. It is clear that ${\hat{Y}}^{-}_{m}=Y^{-}_{m}$. Finally,  we see that ${\hat{Y}}^{+}_{m+2}$ equals $\pm x_g(Y^+_{m+2})$ for $0\le m\le g-2$. 
 \qed 
\enddemo


\bigpagebreak

\flushpar{\bf 5.\quad  Homology generator.}

\medskip

We know from ({0.5}) that
$$\operatorname{H}_{\Cal M}(g-1,e-1,{\pmb \alpha})\cong\operatorname{H}_{\Cal N}(0,0,0)=R.$$Furthermore, the quasi-isomorphism $M\:\Bbb M(eg-e,eg-g)\to \Bbb N(0,0)[-eg]$ is completely given by the map $$M_{g-1}\:\Cal M(g-1,e-1,{\pmb \alpha})\to \Cal N(0,0,0)=R,$$ where 
${M_{g-1}(U \otimes y_1^{(b_1)}\cdots y_g^{(b_g)}\otimes Z)}$ is equal to 
$$\cases [(\omega_E\bowtie y_g^{(e)})\wedge  Z\wedge(U\bowtie (y_1\wedge \dots\wedge y_{g-1}))](\omega_{E^*\otimes G}),
& \text{if $b_g=e-1$, and }\\ 0,& \text{otherwise}. 
\endcases $$ It follows that  $M_{g-1}$
induces an isomorphism $$\operatorname{H}_{\Cal M}(g-1,e-1,\pmb {\pmb \alpha})  \to  R.$$
In this section we identify a cycle $\zeta$ of $\Cal M(g-1,e-1,{\pmb \alpha})$ which is carried to a unit of $R$ by $M_{g-1}$. It follows that the homology class of $\zeta$ generates $\operatorname{H}_{\Cal M}(g-1,e-1,{\pmb \alpha})$. 

\definition {Definition} Let  
$$ \Cal I=\{(i)=(i_1,\dots,i_{e-1})\mid 1\le i_1\le i_2\le \dots\le i_{e-1}\le g\}
$$where each $i_j$ is an integer. 
\flushpar{\bf (a)} Fix $(i)=(i_1,\dots,i_{e-1})\in \Cal I$. For notational convenience, we give meaning to the symbols $i_0$ and $i_{e}$. We allow $i_0$  to mean $1$ and $i_{e}$ to mean $g$. Notice that neither $i_0$ nor $i_{e}$ is an element of the ${e-1}$-tuple $(i)$ . 

\flushpar{\bf (b)} If $(i)\in\Cal I$  and $s$ is an integer, then  $\#_s(i)$ to is  equal to the number of subscripts $p$,   with $1\le p\le e-1$, and $i_p=s$. 

\flushpar{\bf (c)} If $(i)=(i_1,\dots,i_{e-1})$ is in $\Cal I$, then we define $|(i)|$ to be $\sum_{p=1}^{e-1} i_{p}$. 
\enddefinition
\definition {Definition} Let $\zeta$ be the element 
$$  \sum\limits_{(i)\in\Cal I}(-1)^{|(i)|} \prod\limits_{p=1}^{e} u_{p}^{(i_{p}-i_{p-1})}\otimes 
\prod\limits_{s=1}^{g}y_{s}^{(\#_s(i))}
\otimes \left[ \bigwedge\limits_{t=1}^{e}\left(\bigwedge\limits_{w=i_{t-1}}^{i_{t}} (v_{t}\otimes x_{w}) \right)\right] (\omega_{E\otimes G^*}) $$
of $\Cal M(g-1,e-1,{\pmb \alpha})$. 
\enddefinition
\proclaim{Observation}The map $M_{g-1}$ carries $\zeta$ to a unit of $R$.\endproclaim
\demo{Proof}The only term of $\zeta$ which is not killed by $M_{g-1}$ is the term which corresponds to $(i)=(g,\dots,g)$. This term is
$$\pm u_1^{(g-1)}\otimes y_g^{(e-1)}\otimes \left[\left(v_1^{(g)}\bowtie \omega_{G}\right)\wedge\left((v_2\wedge\dots\wedge v_e)\bowtie x_g^{(e-1)}\right)\right](\omega_{E\otimes G^*}),$$
and it is sent to a unit by $M_{g-1}$. \qed
\enddemo

\proclaim{Proposition}The element $\zeta$ of $\Cal M(g-1,e-1,{\pmb \alpha})$ is a cycle and the homology class of $\zeta$ generates  $\operatorname{H}_{\Cal M}(g-1,e-1,{\pmb \alpha})$.
\endproclaim
\demo{Proof} It suffices to show  that $\zeta$ is a cycle. We have $\pmb D(\zeta)$ is equal to $$  \left\{\matrix \format\l\\
\sum\limits_{k,\ell}\sum\limits_{(i)\in\Cal I}(-1)^{|(i)|} v_k\left(\prod\limits_{p=1}^{e} u_{p}^{(i_{p}-i_{p-1})}\right)\otimes 
x_{\ell}\left(\prod\limits_{s=1}^{g}y_{s}^{(\#_s(i))}\right)
\\\hskip25pt \otimes (u_k\otimes y_{\ell})\wedge\left[ \bigwedge\limits_{t=1}^{e}\left(\bigwedge\limits_{w=i_{t-1}}^{i_{t}} (v_{t}\otimes x_{w}) \right)\right] (\omega_{E\otimes G^*}). \endmatrix \right.$$
Notice that the last factor is zero unless $i_{k-1}\le \ell\le i_k$. The middle factor is zero if $i_{k-1}< \ell< i_k$. The first factor is zero if $i_{k-1}=i_k$. For a fixed pair $k,\ell$, the above sum is non-zero provided $i_{k-1}<i_k$ and $\ell$ is equal to one of these two numbers. We have $\pmb D(\zeta)=A+B$, with 
 $$  A=\left\{\matrix \format\l\\
\sum\limits_{k=2}^e\sum\limits_{(i)\in\Cal I\atop{i_{k-1}<i_k}}(-1)^{|(i)|} v_k\left(\prod\limits_{p=1}^{e} u_{p}^{(i_{p}-i_{p-1})}\right)\otimes 
x_{i_{k-1}}\left(\prod\limits_{s=1}^{g}y_{s}^{(\#_s(i))}\right)
\\\hskip25pt \otimes (u_k\otimes y_{i_{k-1}})\wedge\left[ \bigwedge\limits_{t=1}^{e}\left(\bigwedge\limits_{w=i_{t-1}}^{i_{t}} (v_{t}\otimes x_{w}) \right)\right] (\omega_{E\otimes G^*})  \endmatrix \right.$$and $$  B=\left\{\matrix \format\l\\
\sum\limits_{k=1}^{e-1}\sum\limits_{(i)\in\Cal I\atop{i_{k-1}<i_k}}(-1)^{|(i)|} v_k\left(\prod\limits_{p=1}^{e} u_{p}^{(i_{p}-i_{p-1})}\right)\otimes 
x_{i_{k}}\left(\prod\limits_{s=1}^{g}y_{s}^{(\#_s(i))}\right)
\\\hskip25pt \otimes (u_k\otimes y_{i_{k}})\wedge\left[ \bigwedge\limits_{t=1}^{e}\left(\bigwedge\limits_{w=i_{t-1}}^{i_{t}} (v_{t}\otimes x_{w}) \right)\right] (\omega_{E\otimes G^*}). \endmatrix \right.$$
Replace the index $k$ in $A$ by $k+1$. If $k$ is fixed, with $1\le k\le e-1$, then there is a one-to-one correspondence between 
$$\{(i)\in \Cal I\mid i_k<i_{k+1}\}\quad\text{and}\quad  \{(j)\in\Cal I\mid j_{k-1}<j_k\}$$which is given by
$$j_{K}=\cases i_K&\text{if $K\neq k$}\\i_k+1&\text{if $K = k$}.\endcases
$$There is no difficulty in seeing that $A+B=0$. \qed
\enddemo

\bigpagebreak

\flushpar{\bf 6.\quad Homogeneity.}

\medskip

The complexes $\Bbb M(P,Q)$ and $\Bbb N(P,Q)$ each have an enormous amount of homogeneity. This homogeneity passes to the homology modules $\operatorname{H}_{\Cal M}(m,n,p)$ and $\operatorname{H}_{\Cal N}(m,n,p)$, and it even passes across the isomorphism ({0.5}); see Theorem {6.1}. One consequence is Corollary {6.2} which says that 
the homology  modules $\operatorname{H}_{\Cal M}(m,n,p)$ and $\operatorname{H}_{\Cal N}(m,n,p)$ satisfy an extra duality when $e$  is equal to three. This duality is translated in Corollary {6.3} to give an extra symmetry in the graded betti numbers of the  rank one reflexive modules of the determinantal ring defined by the $2\times 2$ minors of a generic $3\times g$ matrix.  

 \definition {Notation}
In the notation of  ({2.1}), let $N=u_1^{\ell_1}\cdots u_e^{\ell_e}$ be a monomial of degree $m+p$ and $M=y_1^{\lambda_1}\cdots y_g^{\lambda_g}$  be a monomial of degree $n+p$. Define $\Cal M(m,n,p)|_{N,M}$ to be  the submodule of $\Cal M(m,n,p)$ which consists of those elements which are homogeneous of degree $\ell_i$ in $u_i$ and of degree $\lambda_j$ in $y_j$ for all $i$ and all $j$. The submodules $\Cal M(m,n,p)|_{N}$ and $\Cal M(m,n,p)|_{M}$, homogeneous  in just the $\{u_i\}$ or just the $\{y_j\}$, are defined in an analogous manner. 
The differential of $\Bbb M$ is homogeneous in the $u$'s and $y$'s; and therefore, the complex $\Bbb M(P,Q)$ naturally decomposes into a direct sum of subcomplexes $\Bbb M(P,Q)|_{N,M}$,  where the sum is taken over all monomials  $N$  and $M$  of degree $P$ and $Q$, respectively. 
Take $\operatorname{H}_{\Cal M}(m,n,p)|_{N,M}$ to mean the homology of the complex $\Bbb M(m+p,n+p)|_{N,M}$ at $\Cal M(m,n,p)|_{N,M}$. We see that $\operatorname{H}_{\Cal M}(m,n,p)$ is equal to $\bigoplus_{N,M}\operatorname{H}_{\Cal M}(m,n,p)|_{N,M}$.  In a similar manner, $\operatorname{H}_{\Cal N}(m,n,p)$ is equal to $\bigoplus_{N',M'}\operatorname{H}_{\Cal N}(m,n,p)|_{N',M'}$,
where the sum is taken over all monomials  $N'$  and $M'$  of degree $m+p$ in the variables $v_1,\dots,v_e$, and degree $n+p$, in the variables $x_1,\dots,x_g$, respectively. \enddefinition

Suppose that the triples $(m,n,p)$ and $(m',n',p')$ satisfy 
$$m+m'=g-1,\quad n+n'=e-1,\quad p+p'={\pmb \alpha},\quad\text{and}\quad  1-e\le m-n\le g-1.$$In this case, ({0.5}) assures us that 
$$\operatorname{H}_{\Cal M}(m,n,p)\cong \operatorname{H}_{\Cal N}(m',n',p').$$
Furthermore, we know that
$$\split \operatorname{H}_{\Cal M}(m,n,p)&=\bigoplus\limits_{N,M}\operatorname{H}_{\Cal M}(m,n,p)|_{N,M}\  \text{and}\\\operatorname{H}_{\Cal N}(m',n',p')&=\bigoplus\limits_{N',M'}\operatorname{H}_{\Cal N}(m',n',p')|_{N',M'},\endsplit$$
as $N$ varies over all monomials of degree $m+p$ in the $u$'s, $M$ varies over all monomials of degree $n+p$ in the $y$'s, $N'$ varies over all monomials of degree $m'+p'$ in the $v$'s, and $M'$ varies over all monomials of degree $n'+p'$ in the $x$'s. It is natural to wonder how a particular pair 
$$\operatorname{H}_{\Cal M}(m,n,p)|_{N,M} \quad\text{and}\quad  \operatorname{H}_{\Cal N}(m',n',p')|_{N',M'}$$ are related.

\proclaim{Theorem {6.1}} Adopt the notation and hypotheses of {\rm({0.5})}. 
Let 
$N=u_1^{a_1}\cdots u_e^{a_e}$, $N'=v_1^{b_1}\cdots v_e^{b_e}$, 
$M=y_1^{c_1}\cdots y_{g}^{c_g}$, and  $M'=x_1^{d_1}\dots x_g^{d_g}$, with $\sum a_i=m+p$, $\sum b_i=m'+p'$, $\sum c_i=n+p$, and $\sum d_i=n'+p'$. If 
$$a_i+b_i=g-1\quad\text{and}\quad  c_j+d_j=e-1$$ for all $i$ and $j$, then 
$$\operatorname{H}_{\Cal M}(m,n,p)|_{N,M} \cong \operatorname{H}_{\Cal N}(m',n',p')|_{N',M'}.$$
\endproclaim
\demo{Proof} We show that the map 
$M_m\:\Cal M(m,n,p)\to\Cal N(m',n',p')$ of Definition {2.6}  carries 
 $$\Cal M(m,n,p)|_{N,M}\quad \text{to} \quad \Cal N(m',n',p')|_{N',M'}.$$ Fix $U\otimes Y\otimes Z\in \Cal M(m,n,p)|_{N,M}$ and  $U''\otimes Y''\otimes Z''\in \Cal M(m',n',p')|_{N'',M''}$, with 
 $$ N''=u_1^{\beta_1}\cdots u_e^{\beta_e}\quad\text{and}\quad   M''=y_1^{\delta_1}\dots y_g^{\delta_g},$$
for  $\sum \beta_i=m'+p'$ and  $\sum \delta_i=n'+p'$.  We know that 
 $[M_m(U\otimes Y\otimes Z)](U''\otimes Y''\otimes Z'')$ is equal to 
$$\left[ \left(\omega_E\bowtie \tau_{m+1}(Y''\otimes Y)\right)\wedge  Z''\wedge  Z \wedge \left(U\bowtie Y^{-}_{m} \right)
 \wedge\left(U'' \bowtie Y^{+}_{m+2}\right) \right](\omega_{E^*\otimes G}).
$$Observe that 
$$\left[ \left(\omega_E\bowtie \tau_{m+1}(Y''\otimes Y)\right)\wedge  Z''\wedge  Z \wedge \left(U\bowtie Y^{-}_{m} \right)
 \wedge\left(U'' \bowtie Y^{+}_{m+2}\right) \right]$$ is an element of the  homogeneous component of ${\tsize \bigwedge}^{eg}(E\otimes G^*)$ which corresponds to the monomials  
$$u_1\cdots u_e\cdot N\cdot N'',\quad y_1\cdot y_g\cdot M\cdot M'',$$and that  $\omega_{E^*\otimes G}$ is an element of the  homogeneous component of ${\tsize \bigwedge}^{eg}(E^*\otimes G)$ which corresponds to the monomials  $$v_1^g\cdots v_e^g,\quad x_1^e\cdots x_g^e.$$ We conclude that if $[M_m(U\otimes Y\otimes Z)](U''\otimes Y''\otimes Z'')$ is not zero, then 
$a_i+\beta_i=g-1$ and $c_j+\delta_j=e-1$ for all $i$ and $j$. \qed
\enddemo

\proclaim{Corollary {6.2}}
{\bf (a)} Assume that $e=3$ and $g$ is arbitrary. If $P$ and $Q$ are integers with $Q-2\le P\le 2Q-1$, then the homology of $\Bbb M(P,Q)$ satisfies 
$$\operatorname{H}_{\Cal M}(m,n,p)\cong\operatorname{H}_{\Cal N}(m',n',p'),$$ provided $$m+m'=Q-1,\quad n+n'=2, \quad\text{and}\quad  p+p'=2Q-2,$$ for $m+p=P$ and $n+p=Q$. 
\flushpar{\bf (b)}
Assume that $g=3$ and $e$ is arbitrary. If $P$ and $Q$ are integers with $P-2\le Q\le 2P-1$, then the homology of $\Bbb M(P,Q)$ satisfies 
$$\operatorname{H}_{\Cal M}(m,n,p)\cong\operatorname{H}_{\Cal N}(m',n',p'),$$ provided $$m+m'=2,\quad n+n'=P-1, \quad\text{and}\quad  p+p'=2P-2,$$ for $m+p=P$ and $n+p=Q$. 
 \endproclaim

\demo{Proof}We prove (a). A symmetric argument establishes  (b).  The $R$-module $\operatorname{H}_{\Cal M}(m,n,p)$ is isomorphic to 
the direct sum of the modules $\operatorname{H}_{\Cal M}(m,n,p)|_{M}$ as $M$ varies over all of the monomials  of degree $Q$ in the variables $y_1,\dots,y_g$. If $M=y_1^{c_1}\cdots y_g^{c_g}$, then Theorem {6.1} tells us that $\operatorname{H}_{\Cal M}(m,n,p)|_{M}$ is zero unless $c_i\le 2$ for all $i$. For each pair of integers $(a,b)$, with $a+2b=Q$, let  $M_{a,b}$ be the monomial $$y_1\cdots y_a\cdot y_{a+1}^2\cdots y_{a+b}^2.$$ If $M$ is any other monomial which consists of $a$ variables raised to the power one and $b$ variables raised to the power $2$, then 
$$\operatorname{H}_{\Cal M}(m,n,p)|_{M} \quad\text{and}\quad  \operatorname{H}_{\Cal M}(m,n,p)|_{M_{a,b}}$$are isomorphic as $R$-modules. Let $\varepsilon_{a,b}$ 
equal the number of monomials in the variables $y_1,\dots,y_g$ which consist of $a$ variables raised to the power one and $b$ variables raised to the power $2$. The exact value of $\varepsilon_{a,b}$ is not needed in this proof; on the other hand, it is obvious that the value is given by 
$$\varepsilon_{a,b}=\cases \frac {g!}{a!b!(g-a-b)!}&\text{if $a+b\le g$}\\ 0&\text{if $g<a+b$}.\endcases$$ We have shown that 
$$\operatorname{H}_{\Cal M}(m,n,p)\cong \bigoplus\limits_{a+2b=Q}\left(\operatorname{Hom}_{\Cal M}(m,n,p)|_{M_{a,b}}\right)^{\varepsilon_{a,b}}.$$ The exact same analysis may be applied to $\operatorname{H}_{\Cal N}(m',n',p')$. We notice that the hypotheses ensure that $n'+p'$ is also equal to $Q$. Let $$M'_{a,b}=x_1\cdots x_a\cdot x_{a+1}^2\cdots x_{a+b}^2.$$It follows that 
$$\operatorname{H}_{\Cal N}(m',n',p')\cong \bigoplus\limits_{a+2b=Q}\left(\operatorname{Hom}_{\Cal N}(m',n',p')|_{M'_{a,b}}\right)^{\varepsilon_{a,b}}.$$
To complete the proof, we apply Theorem {6.1} with $G$ replaced by the free module of rank $Q$ whose basis is $x_1,\dots,x_Q$. It does not matter whether $(x_1,\dots,x_Q)$ is a submodule or an extension of the original $G$. The hypothesis ensures that $-2\le P-Q\le Q-1$. We conclude 
$$\operatorname{Hom}_{\Cal M}(m,n,p)|_{M_{a,b}}\cong \operatorname{Hom}_{\Cal N}(m',n',p')|_{M'},$$ for $M'=x_1\cdots x_a\cdot x_{a+b+1}^2\cdots x_{a+2b}^2$. Notice that for each $i$, with $1\le i\le Q$, the exponent of $y_i$ in $M_{a,b}$ plus the exponent of $x_i$ in $M'$ is equal to $2$.  
 \qed
\enddemo

Take $R$ to be a field (of arbitrary characteristic). We use ({0.2}) to translate Corollary {6.2} into the language of {1.1}. Our base ring is a field, so
$\operatorname{H}_{\Cal N}(m,n,p)\cong\operatorname{H}_{\Cal M}(m,n,p)$.
 If $M$ is a graded $\Cal P$-module,   $\Bbb X\: \dots \to X_1 \to X_0\to M$ is a minimal homogeneous $\Cal P$-free  resolution of $M$,    and $X_p$ is equal to $\bigoplus_i \Cal P(-i)^{\beta_{p,i}(M)}$, then the graded betti number $\beta_{p,q}(M)$ is equal to the dimension of the vector space $\operatorname{Tor}_{p,q}^\Cal P(M,R)$.

\proclaim{Corollary {6.3}} Adopt the data of {1.1} with $R$ equal to a field. 
\flushpar{\bf (a)} Assume that  $e=3$.
If $\ell$ and $q$ are integers with $-2\le \ell\le q-1$, then $$\beta_{p,q}(M_{\ell}) =\beta_{p',q}(M_{\ell'})$$ for $\ell+\ell'=q-3$ and  $p+p'=2q-2$.
\flushpar{\bf (b)} Assume that  $g= 3$.
If $\ell$ and $q$ are integers with $\ell\le 2$ and $1-2\ell \le q$, then  
$$\beta_{p,q}(M_{\ell})=\beta_{p',q'}(M_{\ell'}),$$ provided $$\ell+\ell'=3-\ell-q,\quad q+q'=3(\ell+q-1), \quad\text{and}\quad p+p'=2(\ell+q-1).$$ 
\endproclaim
 
\example{Example} Retain the hypotheses of Corollary {6.3} with $e=g=3$. The graded betti numbers of the $\Cal P$-modules $M_{\ell}$, for $-5\le \ell\le 5$ are
$$\matrix 
^\ddag\beta_{0,5}(M_{-5})=21  &^\ddag\beta_{0,4}(M_{-4})=15 &^\ddag\beta_{0,3}(M_{-3})=10 &\beta_{0,2}(M_{-2})=6\\
^\ddag\beta_{1,6}(M_{-5})=105& ^\ddag\beta_{1,5}(M_{-4})=72 &^\ddag\beta_{1,4}(M_{-3})=45&\beta_{1,3}(M_{-2})=24\\
^\ddag\beta_{2,7}(M_{-5})=216&^\ddag\beta_{2,6}(M_{-4})=141&^\ddag\beta_{2,5}(M_{-3})=81&\beta_{2,4}(M_{-2})=36\\
^\ddag\beta_{3,8}(M_{-5})=234&^\ddag\beta_{3,7}(M_{-4})=144&^\ddag\beta_{3,6}(M_{-3})=74&\beta_{3,5}(M_{-2})=24\\
^\ddag\beta_{4,9}(M_{-5})=141&^\ddag\beta_{4,8}(M_{-4})=81&\beta_{4,7}(M_{-3})=36&\beta_{4,6}(M_{-2})=6\\
^\ddag\beta_{5,10}(M_{-5})=45&\beta_{5,9}(M_{-4})=24&\beta_{5,8}(M_{-3})=9\\
\beta_{6,11}(M_{-5})=6 &\beta_{6,10}(M_{-4})=3&\beta_{6,9}(M_{-3})=1\\
\endmatrix$$
 
$$\matrix 
\beta_{0,1}(M_{-1})=3&\beta_{0,0}(M_{0})=1&\beta_{0,0}(M_{1})=3&\beta_{0,0}(M_{2})=6\\
\beta_{1,2}(M_{-1})=9&\beta_{1,2}(M_{0})=9&\beta_{1,1}(M_{1})=9&\beta_{1,1}(M_{2})=24\\
\beta_{2,3}(M_{-1})=6&\beta_{2,3}(M_{0})=16&\beta_{2,2}(M_{1})=6&\beta_{2,2}(M_{2})=36\\
\beta_{2,4}(M_{-1})=6&\beta_{3,4}(M_{0})=9& \beta_{2,3}(M_{1})=6&\beta_{3,3}(M_{2})=24\\
\beta_{3,5}(M_{-1})=9&\beta_{4,6}(M_{0})=1&\beta_{3,4}(M_{1})=9&\beta_{4,4}(M_{2})=6\\
\beta_{4,6}(M_{-1})=3&&\beta_{4,5}(M_{1})=3\\
\endmatrix$$

$$\matrix 
^\ddag\beta_{0,0}(M_{3})=10&^\ddag\beta_{0,0}(M_{4})=15&^\ddag\beta_{0,0}(M_{5})=21\\
^\ddag\beta_{1,1}(M_{3})=45&^\ddag\beta_{1,1}(M_{4})=72&^\ddag\beta_{1,1}(M_{5})=105\\
^\ddag\beta_{2,2}(M_{3})=81&^\ddag\beta_{2,2}(M_{4})=141&^\ddag\beta_{2,2}(M_{5})=216\\
^\ddag\beta_{3,3}(M_{3})=74&^\ddag\beta_{3,3}(M_{4})=144&^\ddag\beta_{3,3}(M_{5})=234\\
\beta_{4,4}(M_{3})=36&^\ddag\beta_{4,4}(M_{4})=81&^\ddag\beta_{4,4}(M_{5})=141\\
\beta_{5,5}(M_{3})=9&\beta_{5,5}(M_{4})=24&^\ddag\beta_{5,5}(M_{5})=45\\
\beta_{6,6}(M_{3})=1&\beta_{6,6}(M_{4})=3&\beta_{6,6}(M_{5})=6.\\
\endmatrix$$
These numbers were calculated by the computer program Macaulay. There are four symmetries running through this set of numbers. The $R$-modules $E^*$ and $G$ are isomorphic because $e=g$; hence, the $R$-module automorphism of $\Cal P$, which sends the matrix of indeterminates $Z$ to $Z^{\text{\rm T}}$, induces the relation 
$$\beta_{p,q}(M_{\ell})=\beta_{p,q+\ell}(M_{-\ell}),\tag{6.4}$$ for all integers $p$, $q$, and $\ell$. 
If $-2\le \ell\le 2$, then ({1.3}) gives 
$$\beta_{p,q}(M_{\ell})=\beta_{p',q'}(M_{\ell'})\quad\text{for $\ell+\ell'=0$, $p+p'=4$, and $q+q'=6$.}$$
The other two symmetries are listed in Corollary {6.3}. We have marked (with $^\ddag$) the betti numbers  which  satisfy only ({6.4}). Each of the other betti numbers also satisfy at least one of the other symmetries. The module $M_{\ell}$ is Cohen-Macaulay for $-2\le \ell\le 2$, and the symmetries of ({1.3}) apply only in this range. It is interesting to notice that the symmetries of Corollary {6.3} apply to some of the betti numbers of modules $M_{\ell}$ which are not Cohen-Macaulay.
\endexample


\end